\documentclass[12pt]{amsart}
\usepackage{amsbsy}
\usepackage{graphicx,epsfig,subfigure,psfrag}
\begin{document}

\newtheorem{theorem}{Theorem}
\newtheorem{proposition}{Proposition}
\newtheorem{lemma}{Lemma}
\newtheorem{corollary}{Corollary}
\newtheorem{definition}{Definition}
\newtheorem{remark}{Remark}
\newcommand{\beq}{\begin{equation}}
\newcommand{\eeq}{\end{equation}}
\numberwithin{equation}{section} \numberwithin{theorem}{section}
\numberwithin{proposition}{section} \numberwithin{lemma}{section}
\numberwithin{corollary}{section}
\numberwithin{definition}{section} \numberwithin{remark}{section}
\newcommand{\ren}{\mathbb{R}^N}
\newcommand{\re}{\mathbb{R}}
\newcommand{\n}{\nabla}
\newcommand{\iy}{\infty}
\newcommand{\pa}{\partial}
\newcommand{\fp}{\noindent}
\newcommand{\ms}{\medskip\vskip-.1cm}
\newcommand{\mpb}{\medskip}
\newcommand{\BB}{{\bf B}}
\newcommand{\Am}{{\bf A}_{2m}}
\renewcommand{\a}{\alpha}
\renewcommand{\b}{\beta}
\newcommand{\g}{\gamma}
\newcommand{\G}{\Gamma}
\renewcommand{\d}{\delta}
\newcommand{\D}{\Delta}
\newcommand{\e}{\varepsilon}
\newcommand{\var}{\varphi}
\renewcommand{\l}{\lambda}
\renewcommand{\o}{\omega}
\renewcommand{\O}{\Omega}
\newcommand{\s}{\sigma}
\renewcommand{\t}{\tau}
\renewcommand{\th}{\theta}
\newcommand{\z}{\zeta}
\newcommand{\wx}{\widetilde x}
\newcommand{\wt}{\widetilde t}
\newcommand{\noi}{\noindent}
\newcommand{\inA}{\quad \mbox{in} \quad \ren \times \re_+}
\newcommand{\inB}{\quad \mbox{in} \quad}
\newcommand{\inC}{\quad \mbox{in} \quad \re \times \re_+}
\newcommand{\inD}{\quad \mbox{in} \quad \re}
\newcommand{\forA}{\quad \mbox{for} \quad}
\newcommand{\whereA}{,\quad \mbox{where} \quad}
\newcommand{\asA}{\quad \mbox{as} \quad}
\newcommand{\andA}{\quad \mbox{and} \quad}
\newcommand{\ssk}{\smallskip}
\newcommand{\LongA}{\quad \Longrightarrow \quad}
\def\com#1{\fbox{\parbox{6in}{\texttt{#1}}}}

\title
{\bf  Third-order nonlinear dispersion PDEs:\\
 shocks, rarefaction, and blow-up waves}

\author {
Victor A.~Galaktionov and Stanislav I.~Pohozaev}

\address{Department of Mathematical Sciences, University of Bath,
 Bath BA2 7AY, UK}
\email{vag@maths.bath.ac.uk}

\address{Steklov Mathematical Institute,
 Gubkina St. 8, 119991 Moscow, RUSSIA}
\email{pokhozhaev@mi.ras.ru}

\thanks{Research partially supported by the INTAS network
CERN-INTAS00-0136}

  \keywords{Odd-order quasilinear PDE, shock and rarefaction waves,
  entropy solutions,  self-similar patterns.
  {\bf To appear in:} Comp. Math. Math. Phys.
  }
 \subjclass{35K55, 35K65}
 \date{\today}

\begin{abstract}

 Formation of  shocks and collapse of other
discontinuities  via
 smooth rarefaction waves, as well as blow-up phenomena, for the third-order nonlinear
dispersion PDE, as the key model,
 \beq
 \label{01}
 u_t=(uu_x)_{xx} 
 \quad \mbox{in}
 \quad \re \times \re_+,
 \eeq
 are studied.
 Two basic Riemann's problems
  for (\ref{01})
 with initial data
  $$
  S_\mp(x) = \mp{\rm sign}\, x,
  $$
  are shown to create
    the shock ($u(x,t) \equiv S_-(x)$) and smooth rarefaction (for data $S_+$) waves, respectively.
  To this end,
    various blow-up and global similarity solutions of (\ref{01}) revealing the fine structure of
    shock and rarefaction profiles
     are constructed. Next,  eigenfunction and nonlinear capacity techniques for proving blow-up are
  developed.

 The study of (\ref{01}) reveals
 similarities with
 entropy theory for
   scalar conservation laws,
 $$
 u_t + uu_x=0,
  $$
 developed by Oleinik and Kruzhkov
 (equations for $x \in \ren$)
in the 1950-60s.

\end{abstract}

\maketitle


\begin{center}
{\em Dedicated to  the memory of Professors O.A.~Oleinik and
S.N.~Kruzhkov}
\end{center}


\ssk\ssk\ssk

This is the first paper on shock wave theory for higher-order {\em nonlinear dispersion equations}, which will be continued in three more ones \cite{3NDEII, GalNDE5} and \cite{Gal3NDENew}.



\section {Introduction: nonlinear dispersion PDEs and main directions of study}
 \label{Sect1}

 \subsection{NDEs: nonlinear dispersion equations in application and general PDE theory}

 In the present paper,
 we
  develop  basic aspects of singularity and
existence-uniqueness  theory for odd-order {\em nonlinear
dispersion} (or {\em dispersive}) PDEs, which, for short, we call
the NDEs.
The canonical for us  model 
is the 
{\em third-order quadratic NDE} (the NDE--$3$)
 \beq
 \label{1}
  \mbox{$
 u_t={\bf A}(u) \equiv (uu_x)_{xx} = u u_{xxx}+ 3 u_x u_{xx}
 \quad \mbox{in}
 \quad \re \times (0,T), \,\,\, T >0.
  $}
 \eeq
We pose for (\ref{1}) the Cauchy problem
  with locally
integrable initial data
 \beq
 \label{2}
 u(x,0)=u_0(x) \quad \mbox{in} \quad \re.
  \eeq
  Often we assume that $u_0$ is bounded and/or compactly supported.
We will also deal with
the initial-boundary values
problem in $(-L,L)\times \re_+$ with Dirichlet boundary
conditions.

\ssk

\noi\underline{\em On integrable NDEs related to water wave
theory}. Concerning applications of equations such as (\ref{1}),
it is customary that various odd-order PDEs
 appear in classic theory of integrable PDEs,
 such as the {\em KdV equation},
 \beq
 \label{31}
 u_t+uu_x=u_{xxx},
 \eeq
and the {\em fifth-order KdV equation},
 $$
 u_t + u_{xxxxx} + 30 \, u^2 u_x + 20 \, u_x u_{xx} + 10 \, u
 u_{xxx}=0,
 $$
  and
others from shallow water theory.
 The quasilinear
  {\em Harry Dym}
 {\em equation}
  \beq
 \label{HD0}
  u_t = u^3 u_{xxx} \, ,
  \eeq
  which also belongs to the NDE family,
   is
one of the most exotic integrable soliton equations; see
\cite[\S~4.7]{GSVR} for survey and references
 therein.
Integrable equation theory produced various hierarchies of
quasilinear higher-order NDEs, such as
 the fifth-order {\em
Kawamoto equation} \cite{Kaw85}
 $$
 u_t = u^5 u_{xxxxx}+ 5 \,u^4 u_x u_{xxxx}+ 10 \,u^5 u_{xx}
 u_{xxx}.
$$
   Quasilinear integrable extensions are admitted by
 {\em Lax's
 seventh-order KdV equation}
 $$
 \mbox{$
 u_t+\big\{35 u^4 + 70\big[u^2 u_{xx}+ u(u_x)^2\big]+7\big[2 u u_{xxxx}+3 (u_{xx})^2 +
 4 u_x u_{xxx}\big]+u_{xxxxxx}\big\}_x=0,
  $}
  $$
  and by the {\em seventh-order Sawada--Kotara
  equation}
 $$
 \mbox{$
 u_t+\big\{63 u^4 + 63\big[2u^2 u_{xx}+ u(u_x)^2\big]+21\big[ u u_{xxxx}+ (u_{xx})^2 +
  u_x u_{xxx}\big]+u_{xxxxxx}\big\}_x=0;
  $}
  $$
see references in \cite[p.~234]{GSVR}. We will continue this short
survey on existence-uniqueness for some integrable odd-order PDEs
 in Section \ref{Fok1}.

\ssk

\noi\underline{\em Compact pattern formation  and NDEs}. Returning
to lowest third-order NDEs that {\em are not integrable},
 we will also study the   {\em Rosenau--Hyman} (RH)
{\em equation}
  \beq
  \label{Comp.4}
  \mbox{$
  u_t =  (u^2)_{xxx} + (u^2)_x
  $}
  \eeq
  which has special important applications as a widely used model of
   the effects of {nonlinear dispersion} in the pattern
  formation in liquid drops \cite{RosH93}. It is
   the $K(2,2)$ equation from the general $K(m,n)$ family of
   the following NDEs:
  \beq
  \label{Comp.5}
   u_t =  (u^n)_{xxx} +  (u^m)_x \quad (u \ge 0),
   \eeq
   that also
 describe various  phenomena of compact pattern
 formation for $m, \, n>1$, \cite{RosCom94, Ros96}. Such PDEs also appear in curve motion and shortening
flows \cite{Ros00}.
 Similar to well-known parabolic models of porous medium type, the $K(m,n)$ equation (\ref{Comp.5}) with $n>1$
  is degenerated at $u=0$, and
 therefore may exhibit finite speed of propagation and admit solutions with finite
 interfaces.
 The crucial advantage of the RH equation
 (\ref{Comp.4}) is that it possesses
  {\em explicit}  moving compactly supported
 soliton-type solutions, called {\em compactons}
 \cite{RosH93}, which are {\em travelling wave} (TW) solutions.
 We will check entropy properties of  compactons
 for various NDEs.

   Various families of quasilinear
third-order KdV-type equations can be found  in \cite{6R.1}, where
further references concerning such PDEs and their exact solutions
are given.
   Higher-order generalized KdV
equations are of increasing interest; see {e.g.,} the quintic KdV
equation in \cite{Ros98} and \cite{Yao7}, where the seventh-order
PDEs are studied. For the $K(2,2)$ equation (\ref{Comp.4}), the
compacton solutions were constructed in \cite{RosCom94}. 
 More general $B(m,k)$ equations
 (indeed, coinciding with the
 $K(m,k)$
  after scaling)
 $$
 u_t+ a(u^m)_x = \mu (u^k)_{xxx}
 $$
 also admit simple semi-compacton solutions \cite{RosK83}, as well
 as the $Kq(m,\o)$ nonlinear dispersion equation
 (another nonlinear extension of the KdV) \cite{RosCom94}
 $$
 u_t + (u^m)_x + [u^{1-\o}(u^\o u_x)_x]_x =0.
 $$
Setting $m=2$ and $\o=\frac 12$ yields a typical quadratic PDE
 $$
 u_t + (u^2)_x + u u_{xxx} + 2 u_x u_{xx}=0
 $$
 possessing  solutions on standard trigonometric-exponential
   subspaces, where
    $$
    u(x,t)= C_0(t) + C_1(t) \cos \l x + C_2(t)
   \sin \l x
    $$
     and $\{C_0,C_1,C_2\}$ solve a 3D nonlinear dynamical
   system.
 Combining the $K(m,n)$ and $B(m,k)$ equations gives the
 dispersive-dissipativity entity $DD(k,m,n)$ \cite{Ros98DD}
  \index{equation!$DD(k,m,n)$}
  $$
  u_t + a(u^m)_x + (u^n)_{xxx} = \mu (u^k)_{xx}
  $$
  that can also admit solutions on invariant subspaces for
  some values of parameters.

For the fifth-order NDEs,
 such as
  \beq
  \label{Comp.3}
  u_t =  \a (u^2)_{xxxxx} + \b (u^2)_{xxx} + \g (u^2)_x \quad
  \mbox{in} \,\,\, \re \times \re_+,
  \eeq
compacton solutions  were first constructed in \cite{Dey98}, where
the  more general $K(m,n,p)$ family of
PDEs\index{equation!$K(m,n,p)$}
 $$ 
 u_t+ \b_1 (u^m)_x + \b_2 (u^n)_{xxx} + \b_3 D_x^5 (u^p)=0 \quad
 (m,n,p>1)
  $$ 
  was introduced. Some of these equations will be treated later on.
Equation  (\ref{Comp.3}) is also
  associated with the
 family $Q(l,m,n)$ of more general quintic  evolution PDEs with
nonlinear dispersion,\index{equation!$Q(l,m,n)$}
  \beq
  \label{qq1.1}
  \mbox{$
u_t + a (u^{m+1})_{x} + \o \bigl[u(u^n)_{xx}\bigr]_x + \d
\bigl[u(u^l)_{xxxx}\bigr]_x=0,
 $}
 \eeq
 possessing multi-hump, compact solitary solutions \cite{Ros599}.

Concerning higher-order in time quasilinear PDEs, let us mention a
generalization of the {\em combined dissipative double-dispersive}
(CDDD) {\em equation} (see, {e.g.,} \cite{Por02})
 \beq
 \label{Por.1}
 u_{tt}= \a u_{xxxx} + \b u_{xxtt} + 
   \g
 (u^2)_{xxxxt} + \d (u^2)_{xxt} + \e (u^2)_t,
 \eeq
and also the nonlinear modified dispersive Klein--Gordon equation
($mKG(1,n,k)$),
 \beq
 \label{mkg1}
 u_{tt}+ a (u^n)_{xx}+ b(u^k)_{xxxx}=0, \quad n,k>1 \quad (u \ge 0);
  \eeq
 see some exact TW solutions in \cite{Inc07}.
 For $b>0$, (\ref{mkg1}) is of hyperbolic (or Boussinesq) type in
 the class of nonnegative solutions.
Let us also mention a related family of 2D {\em dispersive
Boussinesq equations} denoted by $B(m,n,k,p)$ \cite{Yan03},
 $$ 
  (u^m)_{tt} + \a (u^n)_{xx} + \b (u^k)_{xxxx} + \g
  (u^p)_{yyyy}=0 \quad \mbox{in}
  \quad \re^2 \times \re.
  $$ 
See \cite[Ch.~4-6]{GSVR} for more references and examples of exact
solutions on invariant subspaces of NDEs of various types and
orders.

\ssk

\noi\underline{\em NDEs in general PDE theory}.  In the framework
of general theory of nonlinear evolution  PDEs of the first order
in time, the NDE (\ref{1}) appears the third in the following
ordered list of canonical evolution quasilinear degenerate
equations:
 \beq
 \label{ca1}
  \mbox{$
 u_t= - \frac 12\, (u^2)_x \quad (\mbox{the conservation law}),
 \qquad\qquad\qquad\qquad\qquad\,\,\,\,
  $}
  \eeq
 \beq
 \label{ca2}
  \mbox{$
 u_t= \frac 12 \,(u^2)_{xx}\quad (\mbox{the porous medium equation}),
 \qquad\qquad\qquad\quad\,\,\,\,
 $}
 \eeq
 \beq
\label{ca3} \mbox{$
 u_t= \frac 12 \,(u^2)_{xxx}\quad (\mbox{the nonlinear dispersion equation}),
 \qquad\qquad\,\,\,\,\,\,\,
 $}
 \eeq
  \beq
\label{ca4}
 \mbox{$
 u_t= -\frac 12 \,(|u|u)_{xxxx}\quad (\mbox{the $4^{\rm th}$-order nonlinear diffusion
 equation}).
 $}
 \eeq
 In (\ref{ca4}), the quadratic nonlinearity $u^2$ is replaced by
 the monotone one $|u|u$ in order to keep the parabolicity on
 solutions of changing sign. The same can be done in the PME
 (\ref{ca2}), though this classic version is parabolic on
 nonnegative solutions, a property that is preserved by the
 Maximum Principle.
The further extension of the list by including
  \beq
   \label{ca5}
    \mbox{$
 u_t= -\frac 12 \,(u^2)_{xxxxx}\quad (\mbox{the NDE--5}) \quad \mbox{and}
 \qquad\qquad\qquad\qquad\qquad\,\,\,\,
 $}
 \eeq
 \beq
\label{ca6}
 \mbox{$
 u_t= \frac 12 \,(|u|u)_{xxxxxx}\quad (\mbox{the $6^{\rm th}$-order nonlinear diffusion
 equation})
 $}
 \eeq
is not that essential since these PDEs belong to the same families
as (\ref{ca3}) and (\ref{ca4}) respectively with similar covering
mathematical concepts (but indeed more difficult).

\ssk

 Mathematical theory of first two equations, (\ref{ca1}) (see
detailed survey and references below) and (\ref{ca2}) (for quoting
 main achievements of  PME theory developed in the 1950--80s,
see e.g., \cite[Ch.~2]{AMGV}), was essentially completed in the
twentieth century. It is curious that looking more difficult the
fourth-order nonlinear diffusion equation (\ref{ca4}) has a
monotone operator in $H^{-2}$, so the Cauchy problem admits a
unique weak solution as follows from classic theory; see Lions
\cite[Ch.~2]{LIO}. Of course, some other qualitative properties of
solutions of (\ref{ca4}) are  more difficult and remain open
still.

It turns out that, rather surprisingly,  the third order NDE
(\ref{ca3}) is the only one in this basic list that has rather
obscure understanding and lacking of a reliable mathematical basis
concerning generic singularities, shocks, rarefaction waves, and
entropy-like theory.

\subsection{Mathematical preliminaries: analogies with conservation laws, Riemann's problems for
basic shocks $S_\mp(x)$ and $H(\mp x)$, and first results}
 \label{Sect1.2}

  As a key feature of our analysis,
 equation
(\ref{1}) inherits clear similarities 
of the behaviour for the first-order conservation laws such as
{\em Euler's equation} (same as (\ref{ca1})) from gas dynamics
 \beq
 \label{3}
  u_t + uu_x=0
\quad \mbox{in} \quad \re\times \re_+,
 \eeq
 whose entropy theory  was created by Oleinik \cite{Ol1, Ol59} and Kruzhkov
\cite{Kru2} (equations in $\ren$) in the 1950--60s; see details on
the history, main results, and modern developments in the
well-known monographs \cite{Bres, Daf, Sm}\footnote{First results
on the formation of shocks were performed by Riemann in 1858
\cite{Ri58}; see \cite{Chr07}.}.

\ssk

As for (\ref{3}), in view of the full divergence of the equation
(\ref{1}), it is natural to define weak solutions. For
convenience, we present here  a standard definition mentioning
that, in fact, the concept of weak solutions for NDEs even in
fully divergent form is not entirely consistent, to say nothing
about other non-divergent equations admitting no weak formulation
at all.


\ssk

\noi{\bf Definition \ref{Sect1}.1.} A function $u=u(x,t)$ is a
weak solution of $(\ref{1})$, $(\ref{2})$ if 

 (i) $u^2 \in L^1_{\rm loc}(\re\times (0,T))$,

 (ii) $u$ satisfies (\ref{1})
in the weak sense:
for any test function $\varphi(x,t) \in C_0^\infty(\re \times
(0,T))$,
 \beq
 \label{21}
 \mbox{$
 \iint u \varphi_t=\frac 12\,  \iint u^2 \varphi_{xxx},
 $}
  \eeq
and (iii) satisfies the initial condition (\ref{2}) in the sense
of distributions,
 \beq
 \label{22}
 \mbox{$
{\rm ess} \,  \lim_{t \to 0}\int u(x,t) \psi(x) = \int u_0(x) \psi
(x) \quad
 \mbox{for any} \quad \psi \in C_0^\infty(\re).
  $}
   \eeq

   We will show that, for some data, (\ref{2}) is also true
    with convergence in $L^1_{\rm loc}(\re)$. This topology is
natural for conservation laws (\ref{3}), but for the NDE (\ref{1})
 such a convergence is not straightforward and demands a precise
 knowing of the structure of oscillatory ``tails" of similarity
 rarefaction solutions to be studied in detail.
The assumption $T<\infty$ is often essential, since, unlike
(\ref{3}), the NDE (\ref{1}) can produce complete blow-up from
bounded data.

\ssk

Thus, again similar to (\ref{3}), one observes a typical
difficulty: according to Definition \ref{Sect1}.1,
 both discontinuous step-like functions
  \beq
  \label{Ri1}
  S_\mp(x) = \mp{\rm sign}\, x = \mp \left\{
 \begin{matrix}\,\,\,
 1 \,\,\, \mbox{for} \,\,\, x>0,\\
\,\,\,  0 \,\,\, \mbox{for} \,\,\, x=0,\\ -1
 \,\, \mbox{for}
\,\,\, x < 0,
 \end{matrix}
 \right.
  \eeq
 are weak {\em stationary}  solutions of (\ref{1}) satisfying
 \beq
 \label{Ri33}
 (u^2)_{xxx}=0
  \eeq
  in the  weak sense,
 since in (\ref{21}) $u^2(x)=S_\mp^2(x) \equiv 1$ is $C^3$
smooth (and analytic).

Again referring to entropy theory  for conservation laws
(\ref{3}),  it is well-known that
 \beq
 \label{rr1}
 \begin{matrix}
 S_-(x) \,\,\, \mbox{is the entropy shock wave, and}\ssk\\
 S_+(x) \,\,\, \mbox{is not an entropy solution}.\qquad
 \end{matrix}
  \eeq
This means that
 \beq
 \label{rr41}
 u_-(x,t) \equiv S_-(x)=-{\rm sign} \, x
  \eeq
is the unique entropy solution of the PDE (\ref{3}) with the same
initial data $S_-(x)$. On the contrary, taking $S_+$ initial data
  yields the   {\em rarefaction wave}
    with a simple similarity piece-wise linear
    structure
   \beq
   \label{rr6}
   \mbox{$
    u_0(x)=S_+(x)={\rm sign} \, x \,\,\, \Longrightarrow \,\,\,
   u_{+}(x,t)= g(\frac xt) = \left\{ \begin{matrix}
    -1 \quad \mbox{for} \,\,\, x<-t, \cr
 \,\,\frac xt \quad  \mbox{for} \,\,\, |x|<t, \cr
  1 \quad \,  \mbox{for} \,\,\, x>t.
 \end{matrix}
   \right.
    $}
 \eeq

\ssk

\noi\underline{\em First aims of singularity theory}. Thus, as in
classic conservation law theory, the above discussion formulate
our first two aims:

\ssk

$\bullet$ Our first one is to justify the same classification of
main two Riemann's problems with data (\ref{Ri1}) for the NDE
(\ref{1}) and to construct the corresponding rarefaction wave for
$S_+(x)$, as an analogy of (\ref{rr6}) for the conservation law;

\ssk

 $\bullet$ The second aim is to describe evolution blow-up
formation for (\ref{1})
 of another shock that is the {\em reflected Heaviside function} given
 by
 \beq
 \label{2.12}
 H(-x)= \left\{
  \begin{matrix}
  1 \,\,\,\, \mbox{for} \,\,\,\, x<0,\ssk \\
 0 \,\,\,\, \mbox{for} \,\,\,\, x>0\,\,
  \end{matrix}
   \right.
    \eeq
    ($H(0)=0.4197...$ will be determined by construction).
    It turns out that this third Riemann problem for (\ref{1})
is also solved in a blow-up similarity fashion with the solution
exhibiting a {\em finite interface}. As a crucial feature, note
that $H(-x)$ {\em is not} a weak stationary solution of (\ref{1}).
In general, we therefore would like to claim the following rule:
 \beq
 \label{RW}
 \mbox{shock waves of higher-order NDEs are not necessarily weak
 solutions.}
  \eeq
Of course, this is obvious for non-divergent NDEs or those that
are fully nonlinear (we will consider such equations), but this
also happens for fully divergent PDEs such as (\ref{1}).


\subsection{NDEs are not hyperbolic systems}

 It should be noted that by no means (\ref{1})
can be treated as a hyperbolic system, for which   fully developed
theory in 1D is now available; see e.g., Bressan \cite{Bres}.
Indeed, rewriting (\ref{1}) as a first-order system yields
 \beq
 \label{4}
 \left\{
 \begin{matrix}
u_t=w_x, \\
   v= uu_x,\\
    w=v_x.\,\,
    \end{matrix}
    \right.
     \eeq
Evidently, it is not a hyperbolic system that represents a
combination of a single first-order evolution PDE with two
stationary equations. In a natural formal sense,
(\ref{4}) can be considered  as the limit as $\e \to 0$ of a more
standard first-order system such as
 \beq
 \label{5}
 \left\{
 \begin{matrix}
u_t=w_x, \qquad\quad\,\,\,\ssk\\
   v_t= \frac 1\e\,(v-u u_x),\ssk\\
    w_t= \frac 1\e \,(w-v_x).
    \end{matrix}
    \right.
     \eeq
Nevertheless, the principal differential operator in (\ref{5}) is
not hyperbolic for any $|\e| \ll 1$ that is easily seen by
calculating  eigenvalues of the matrix,
 \beq
 \label{6}
 \left|
 \begin{matrix}
 -\l \,\,\,\,\,\,0\, \,\,\,\,\,\,\, 1 \\
   -\frac u\e \, -\l \,\,\,\,\,0\\
    \,\,0 \,\, - \frac 1\e \, -\l
    \end{matrix}
    \right|=0
    \quad
    \mbox{$
     \Longrightarrow \quad \l^3= \frac u{\e^2}.
     $}
     \eeq
The algebraic equation in (\ref{6})  admits complex roots for any
$u \not =0$ that reflects the actual highly oscillatory properties
of solutions of both systems (\ref{6}) and (\ref{5}) and the
original PDE (\ref{1}). Note that the generic behaviour of small
solutions of the KdV equation (\ref{31}) with linear dispersion
mechanism is also oscillatory. This is associated with  the
fundamental solution of the corresponding operator $D_t-D_x^3$
given by Airy function of changing sign, which is not absolutely
integrable on the whole line $x \in\re$ for any $t>0$.

\subsection{Conservation law in $H^{-1}$}

On the other hand, equation (\ref{1}) can be written as
 \beq
 \label{10}
 (-D_x^2)^{-1}u_t + uu_x=0,
  \eeq
 with the following usual definition of the linear operator
 $(-D_x^2)^{-1}>0$ in $\re$:
 \beq
 \label{11}
 \begin{matrix}
 (-D_x^2)^{-1} v=g, \quad \mbox{if} \quad g''=-v, \,\,\,g(\pm
 \infty)=0 \quad \big( \int\limits_{-\infty}^{+\infty}v=0\big), \quad \mbox{i.e.,}\qquad\ssk\ssk\ssk\\
   g(x)= \int\limits_x^\infty \int\limits_{-\infty}^y v(z) \, {\mathrm d}z\, {\mathrm
   d}y, \quad \mbox{provided that}
   \quad v \in L^1(\re), \,\,\, \int\limits_{-\infty}^x v \in
   L^1(\re).\qquad
    \end{matrix}
     \eeq
In a bounded interval $[-L,L]$ to be also used, the definition is
standard, with $g(\pm L)=0$.

The form (\ref{10}) makes it possible to get the first {\em a
priori} uniform bound on solutions for data $u_0 \in H^{-1}(\re)$:
multiplying (\ref{10}) by $u$ in $L^2$ we get the conservation law
 \beq
 \label{111}
  \mbox{$
 \frac 12 \, \frac{\mathrm d}{{\mathrm d}t}\, \|u(t)\|_{H^{-1}}^2
 =0
 \quad \Longrightarrow \quad \|u(t)\|_{H^{-1}}^2= c_0= \|u_0\|_{H^{-1}}^2
 \quad \mbox{for all} \,\,\, t>0.
  $}
  \eeq
Notice that this is weaker than the estimate for the conservation
law (\ref{3}), for which (\ref{111}) takes place in the topology
of $L^2(\re)$.

It seems that the representation (\ref{10}) and the estimate
(\ref{111}) do not help neither better understanding of formation
of singularities nor entropy essence of solutions.
 Nevertheless, (\ref{10}) convinces to classify the NDE
(\ref{1}) as  a {\em non-local conservation law}  posed  in
$H^{-1}$.

\subsection{Plan of the paper:
 known results and
  main directions of our singularity study}
 \label{Fok1}

We mention again that, surprisingly,  quite a little is known
about suitable mathematics of shocks, rarefaction waves, and
entropy essence of the NDEs such as (\ref{1}).
Moreover, even for sufficiently smooth, continuous solutions, it
is not exaggeration to say that, for (\ref{1}) and the compacton
equation (\ref{Comp.4}), even basic facts concerning proper posing
the Cauchy problem remained obscure, especially, existence of
discontinuous shock\footnote{
   It is relevant commenting that a few years ago,
   one of the first author's paper on shock waves for the NDE (\ref{1})
 submitted to a respectable journal partially dealing with nonlinear PDEs
 received
   the accusation
   that ``... this equation
 cannot admit discontinuous
 solutions at all as the heat equation". It seems the Referee was
 confused observing in (\ref{1}) too many (three) spatial
 derivatives, while the heat equation $u_t=u_{xx}$ has two only.
Of course,  any odd, $2m+1$ (say, 2007 for $m=1003$), number of
derivatives
 in $u_t=(-1)^{m+1}D_x^{2m} (u u_x)$
  will lead to discontinuous
shocks for some initial data.}  and rarefaction
   waves, as well as {\em compactons} as compactly
supported solutions (in the twenty-first century!).

In the present paper, we describe various types of shock,
rarefaction, and blow-up waves for NDEs under consideration.
 In the next paper \cite{3NDEII},
 we are going to propose key concepts for developing
adequate mathematics of NDEs with shocks, which will be concluded
 by revealing connections with other classes
of nonlinear degenerate PDEs. It turns out that some NDE concepts
has definite reliable common roots and can be put into the
framework of much better developed theory of quasilinear parabolic
equations.

\ssk

\noi\underline{\em Back to integrable models: existence and
uniqueness}. In fact, modern  mathematical theory of odd-order
quasilinear PDEs is originated and continues to be strongly
connected with the class of integrable equations. Special
advantages of integrability by using the
inverse scattering transform method, Lax pairs, Liouville
transformations,
 and other explicit algebraic manipulations made it
possible to create rather complete theory for some of these
difficult quasilinear PDEs. Nowadays,  well-developed theory
 and  most of rigorous results on existence,
uniqueness, and various singularity and non-differentiability
properties are associated with NDE-type integrable models such as
{\em Fuchssteiner--Fokas--Camassa--Holm} (FFCH) {\em
 equation}
 \beq
  \label{R1}
  \mbox{$
  (I-D_x^2)u_t=
 - 3u u_x+ 2u_x u_{xx} + u u_{xxx}
  \equiv -(I-D_x^2)(u u_x)- \big[u^2 + \frac 12\, (u_x)^2\big].
  $}
\eeq
  Equation (\ref{R1})
 is   an asymptotic model  describing  the wave dynamics
at the free surface of fluids under gravity. It is derived from
Euler equations for inviscid fluids under the long wave
asymptotics of shallow water behavior (where the function $u$ is
the height of the water above a flat bottom).
 Applying to (\ref{R1}) the integral operator $(I-D_x^2)^{-1}$
 with the $L^2$-kernel $\o(s)= \frac 12 \, {\mathrm
 e}^{-|s|}>0$, reduces it, for a class of solutions, to the
conservation law (\ref{1}) with a compact {\em first-order}
perturbation,
 \beq
 \label{Con.eq}
  \mbox{$
 u_t  + u u_x=  - \big[ \o*\big(u^2 + \frac 12  (u_x)^2 \big)
 \big]_x.
  $}
 \eeq
 Almost all mathematical results (including entropy
 inequalities and Oleinik's condition (E)) have been obtained by
 using this integral representation of the FFCM equation;
 see a long list of references given in \cite[p.~232]{GSVR}.

There is  another  integrable PDE from the family with third-order
quadratic operators,
 \beq
 \label{ff1SSH}
 u_t-u_{xxt}= \a u u_x + \b u_x u_{xx} + u u_{xxx} \quad (\a, \,
 \b \in \re),
  \eeq
  where $\a=-3$ and $\b=2$ yields the FFCH equation (\ref{R1}).
This is  the {\em Degasperis--Procesi equation} for $\a=-4$ and
$\b=3$,
 \beq
 \label{DP0}
 u_t - u_{xxt}= -4 u u_x+ 3u_x u_{xx} + u u_{xxx}.
 \eeq
On existence, uniqueness (of entropy solutions in $L^1 \cap BV$),
parabolic $\e$-regularization, Oleinik's entropy estimate, and
generalized PDEs, see \cite{Coc06}.
 Besides (\ref{R1}) and (\ref{DP0}), the family
(\ref{ff1SSH}) does not contain other integrable entries.
 A  list of
more applied papers related to various NDEs is also available in
\cite[Ch.~4]{GSVR}.

\ssk

\noi\underline{\em Main directions of study}.
 Concerning the simple canonical  model (\ref{1}), which is not integrable
 and does not admit a reduction
 like (\ref{Con.eq}) with a compact first-order perturbation,
 we  do
 the following research:

 \ssk

$\bullet$ {\bf (i)} To check whether (\ref{1}) admits
discontinuous {\em shock
 wave}
 solutions, and which smooth  similarity solutions can create
 shocks in finite time.

 \ssk

$\bullet$  {\bf (ii)} To verify which discontinuous solutions are
{\em entropy} in a natural sense,
  and which are not that give rise to smooth {\em rarefaction
 waves}. To this end, we apply the idea of smooth deformations of
 shocks to see which ones are stable.
    In particular,
 we show  that
 two basic Riemann's problems
  for (\ref{1})
 with initial data (\ref{Ri1})
   correspond
    the shock ($S_-$) and  rarefaction ($S_+$) waves respectively.
This coincide with the classification for the conservation law
(\ref{3}). We also solve the third Riemann problem with the
Heaviside function (\ref{2.12}) that turns out to be
entropy. Similarly, the reflected data $H(x)$ are shown to produce
 the corresponding smooth rarefaction wave.

\ssk

 $\bullet$ {\bf (iii)} For these purposes, construct various self-similar solutions
 of (\ref{1}) describing formation of shocks and collapse of
 non-evolution discontinuities (i.e., leading to rarefaction waves); here we
 also
 use Gel'fand's concept (1959) \cite[\S\S~2,\,8]{Gel} of G-admissible solutions of
 higher-order ODEs.

 \ssk

 $\bullet$  {\bf (iv)} As a consequence, to prove that
  (\ref{1}) describes processes with finite propagation of
  perturbations (indeed, this was well-known for decades but, it seems,  was
not suitably treated mathematically).

\ssk

In the forthcoming paper \cite{3NDEII}, we continue our study and:

\ssk

 $\bullet$ {\bf (v)} Develop local existence and uniqueness theory for
 the NDE (\ref{1}) and introduce the concept of the $\d$-entropy test.

\ssk

$\bullet$  {\bf (vi)}  For the RH equation such as (\ref{Comp.4}),
we prove that Rosenau's compacton solutions are both entropy and
G-admissible.

\ssk

 Some of related  questions and results were previously discussed in
 a more applied and formal fashion in \cite[\S~7]{GalEng} and
 \cite[Ch.~4]{GSVR}, where another idea of
 regular approximations was under scrutiny. Namely,
 {\em proper solutions} therein were obtained in the limit
  $\e \to 0^+$ of classical
 solutions $\{u_\e\}$ of a family of uniformly parabolic
equations with same data
 \beq
 \label{ee1}
  u_t=(uu_x)_{xx}- \e u_{xxxx} \quad \mbox{in}
 \quad \re \times \re_+ \quad (\e>0).
  \eeq
In particular, it was shown that a direct verification that the
$\e$-approximation (\ref{ee1}) yields as $\e \to 0$ the correct
Kruzhkov's-type entropy solution leads to difficult open problems.

Further extensions to other three- and higher-order NDEs can be
found in \cite{GalNDE5, 3NDEII}.

\section{First blow-up results by two methods}
 \label{SectBl1}

In this section, we begin our study of blow-up singularities that
can be generated by the NDE (\ref{1}) and others. We pose an IBVP
and demonstrate two different methods to prove global nonexistence
of a class of solutions. Note that these results differ from those
that have been known in the theory of more classic {\em semilinear
dispersive equations} (SDEs)
 such as
 \beq
 \label{sde1}
 \mbox{$
  {\rm i} \, u_t + R(-{\rm i}\, D_x)u - Q(-{\rm i}\,
  D_x)(u^{p+1})=0, \quad p \ge 1,
   $}
   \eeq
   where, typically, the functions $P$ and $P$ are polynomials in
   one variables with real coefficients;
 see \cite{BonaW1} for the results and
a survey.

 \subsection{A modified eigenfunction method}

 We begin by deriving simple ordinary
differential inequalities (ODIs) for Fourier-like expansion
coefficients of solutions. Here we perform a convenient adaptation
of Kaplan's eigenfunction method (1963) to odd, $(2m+1)$th-order
nonlinear PDEs. This completes our analysis performed for  even,
$2m$th-order parabolic equations
 in \cite[\S~5]{GP1}, where key references and other
 basic
results can be found.
 Another alternative approach to nonexistence for such PDEs is
 developed in author's resent paper \cite{GPNon}.
 Later on, in Section \ref{Sect3}, we will
 show that looking analogous $L^2$-blow-up is actually possible and can
 be
 driven by self-similar solutions.

For convenience and simplicity,  we  consider (\ref{1}) for $x<0$
and  take initial data
 \beq
 \label{991}
 u_0(x) \le 0 \quad (u_0(x) \not \equiv 0).
 \eeq
We also assume the following conditions at the origin:
 \beq
  \label{992}
u=0, \,\,\,(u^2)_x \ge 0 \quad \mbox{and the flux} \,\,\,
(u^2)_{xx} \le 0 \,\,\,
  \mbox{at} \,\,\, x=0.
   \eeq
Note that all these conditions hold for the Cauchy problem if the
right-hand interface of the given compactly supported solution
belong to $\{x<0\}$ for $t \in (0,T)$. We will describe such
solutions with finite interfaces later on. On the other hand, we
may pose an initial-boundary value problem (IBVP) for (\ref{1}) in
$\re_- \times (0,T)$ with the necessary conditions at $x=0$ such
that (\ref{992}) holds. Not discussing here the local
 well-posedness of such an IBVP, we will show that some of its solutions
  will then blow-up in finite time.

We next fix a constant $L>0$ and choose the ``eigenfunction" as
follows:
 \beq
 \label{993}
 \phi(x)=-(x+L)^3 < 0 \,\,\, \mbox{for} \,\,\ x \in (0,L),
 \quad
 \phi(-L)=\phi'(-L)=\phi''(-L)=0.
   \eeq
Here $\phi$ satisfies the obvious ``eigenfunction" equation
 \beq
 \label{777}
{\bf P}_3 \phi \equiv \phi'''(x) = -6 \quad \mbox{on} \quad (-L,0)
\quad (\phi \le 0),
 \eeq
 plus the above boundary conditions. Here the linear operator
 ${\bf P}_3=D_x^3$ is naturally associated with the NDE (\ref{1}).
Actually, any $\phi$ satisfying the differential inequality
 \beq
 \label{778}
\phi'''(x) \le -C_0 \quad (\phi \le 0),
 \eeq
 where $C_0>0$ is a constant,
would successfully fit into our further analysis.

   Let us introduce the corresponding expansion coefficient
  of the solution
   \beq
   \label{994}
    \mbox{$
    \begin{matrix}
   J(t) = \int\limits_{-L}^0 u(x,t) \phi(x)\,{\mathrm d}x \equiv
-\int\limits_{-L}^0 u(x,t) (x+L)^3\,{\mathrm d}x, \ssk \\
\mbox{where} \quad J(0)=-\int\limits_{-L}^0 u_0(x)
(x+L)^3\,{\mathrm d}x>0.
 \end{matrix}
 $}
 \eeq
Multiplying (\ref{1}) by $\phi(x)$ and integrating over $(-L,0)$
by parts three times, by using the properties in (\ref{992}),
(\ref{993}), and (\ref{777}), one obtains that
 \beq
 \label{995}
  \mbox{$
  \begin{matrix}
 J'(t) = \frac 12\big[(u^2)_{xx} \phi-(u^2)_{x} \phi_x+
  u^2\, \phi_{xx}\big]_{-L}^0- \frac 12 \, \int\limits_{-L}^0 u^2(x,t)\phi'''(x)\,{\mathrm d}x
  \ssk \\
   \ge - \frac 12 \, \int\limits_{-L}^0 u^2(x,t)\phi'''(x)\,{\mathrm d}x
  =3\, \int\limits_{-L}^0 u^2(x,t)\,{\mathrm d}x.
 \end{matrix}
  $}
  \eeq
We next use the H\"older inequality to estimate the right-hand
side in (\ref{995}):
 $$
   \begin{matrix}
3\, \int\limits_{-L}^0 u^2(x,t)\,{\mathrm d}x \ge
   3\, \int\limits_{-L}^0
u^2(x,t) \frac{|x+L|^6}{L^6}\,{\mathrm d}x
 \ssk \\
  =
 \frac 3{L^5}\, \int\limits_{-L}^0 (|x+L|^3 |u(x,t)|)^2 \frac 1 L\,{\mathrm d}x
  \ge
  \frac 3{L^7}\, \big(\int\limits_{-L}^0 |x+L|^3 |u(x,t)| \,{\mathrm d}x
  \big)^2 \ge \frac 3{L^7} J^2(t).
   \end{matrix}
  $$

   Thus,
   we obtain the
   following simple quadratic ODI for the expansion coefficient $J(t)$:
 \beq
 \label{997}
  \mbox{$
 J'(t) \ge \frac 3{L^7} J^2(t) \,\,\, \mbox{for} \,\,\, t>0, \quad
 J(0)=J_0>0.
  $}
  \eeq
Integration over $(0,t)$ yields a typical blow-up,
 \beq
 \label{998}
  \mbox{$
  J(t) \ge \frac{L^3}{3(T_0-t)} \quad \mbox{with the blow-up time}
   \quad T \le T_0= \frac{L^7}{3J_0}.
    $}
     \eeq
Recall again that this blow-up can happen if the flux $(u^2)_{xx}$
at $x=0$ satisfies the sign condition in (\ref{992}). Otherwise,
if it is violated and the flux gets positive, this can prevent any
blow-up in $L^2$ (but cannot prevent a ``weaker" gradient blow-up
of $u_x$ to be studied in greater detail).

This approach to blow-up admits natural adaptation to fifth- and
higher-order NDEs
 \beq
 \label{3,1}
 u_t=(u u_{xx})_x \equiv u u_{xxx} + u_x u_{xx},
  \eeq
  which we call the NDE--(2,1), where 2 and 1 stand for the number
  of the internal and external derivatives in this differential
  form.
 For second- and third-order in time NDEs
\beq
   \label{222}
   u_{tt}=(uu_x)_{xx} \quad \mbox{and}
   \quad  u_{ttt}=(uu_x)_{xx}.
    \eeq
 the same
integration by parts leads, instead of (\ref{997}), to ODIs
 \beq
 \label{9981}
 \mbox{$
 J''(t) \ge \frac 3{L^7} J^2(t) \quad \mbox{and}
 \quad
J'''(t) \ge \frac 3{L^7} J^2(t).
  $}
  \eeq
Here blow-up depends on initial values $J_0=J(0)$, $J'(0)$, and
$J''(0)$. Finally,
  these formally define respectively the following  blow-up rates:
 \beq
 \label{999}
 \mbox{$
 J(t) \sim \frac 1{(T_0-t)^2} \quad \mbox{and}
 \quad  J(t) \sim \frac 1{(T_0-t)^3}.
  $}
   \eeq

\subsection{Nonlinear capacity method: another version of blow-up}

We consider the second-order in time NDE in (\ref{222}) in $\re_+
\times \re_+$ ,
 \beq
 \label{881}
 \mbox{$
 u_{tt}= \frac 12\,(u^2)_{xxx}
 $}
 \eeq
and, for simplicity, assume all three {\em zero conditions}
(\ref{992}) at the origin. We now choose a different cut-off
function,
 $$
 \phi(x,t)=\phi_0(t)(L-x)^3, \,\,\, \mbox{where}
  \,\,\, \phi_0 \ge 0, \,\, \phi_0(0)=\phi_0'(0)=0,\,\,
  \phi_0(T)=\phi_0'(T)=0.
  $$
Integrating by parts yields
 $$
 \begin{matrix}
  \frac 12\int\limits_0^L (u^2)_{xxx} \phi(x,t) \, {\mathrm d}x= 3 \int\limits_0^L u^2
  \phi_0(t)\, {\mathrm d}x, \\
  \int\limits_0^T \int\limits_0^L u_{tt}\phi\,{\mathrm d}x\,{\mathrm d}t=
  - \int\limits_0^L u_t(x,0)(x-L)^3\, {\mathrm d}x+
  \int\limits_0^T \int\limits_0^L u \phi_0''(t)(x-L)^3\,{\mathrm d}x\,{\mathrm
  d}t.
  \end{matrix}
  $$
 Hence, we obtain the integral identity
\beq
 \label{882}
 \mbox{$
3 \int\limits_0^T \int\limits_0^L u^2
  \phi_0(t)\, {\mathrm d}x\, {\mathrm d}t= - \int\limits_0^L u_t(x,0)(x-L)^3
  \, {\mathrm d}x + \int\limits_0^T \int\limits_0^L u \phi_0''(t)(x-L)^3\,{\mathrm d}x\,{\mathrm
  d}t.
 $}
 \eeq
By Young's inequality, setting for convenience
$\phi_1(x)=(x-L)^3$,
 $$
 \mbox{$
 \iint u \phi_0''\phi_1 = \iint u \sqrt{\phi_0}\,\,
 \frac{\phi_0''\phi_1}{\sqrt{\phi_0}} \le
   \frac 12 \iint u^2 \phi_0 + \frac 12 \iint
   \frac{|\phi_0''|^2}{\phi_0} \phi_1^2,
 $}
   $$
   so that (\ref{882}) implies
 \beq
 \label{885}
  \mbox{$
   \frac 52 \int\limits_0^T \int\limits_0^L u^2 \phi_0 \le \frac 12 \int\limits_0^T \int\limits_0^L \frac{|\phi_0''|^2}{\phi_0}
   \phi_1^2-
    \int\limits_0^L u_t(x,0) \phi_1(x) \, {\mathrm d}x.
      $}
 \eeq
 We next perform scaling
  $$
   \begin{matrix}
  \mbox{$
   t \mapsto \t=\frac t T \quad \mbox{and replace} \quad
 \phi_0(t) \mapsto \tilde \phi_0(\t), \quad \mbox{and then}
  $} \ssk\\
 \int\limits_0^T \int\limits_0^L \frac{|\phi_0''(t)|^2}{\phi_0(t)} \phi_1^2(x)\, {\mathrm d}x {\mathrm
 d}t
  = \int\limits_0^T \frac{|\phi_0''(t)|^2}{\phi_0(t)}\, {\mathrm d}t \int\limits_0^L \phi_1(x)\, {\mathrm d}x
  \ssk\\
  = \frac{L^7}7\int\limits_0^T \frac{|\phi_0''(t)|^2}{\phi_0(t)}\, {\mathrm
  d}t\big|_{t \mapsto \t=\frac tT} = \frac{L^7}7 \, \frac 1{T^3} c_0,
 \,\,\,  c_0= \int\limits_0^1 \frac{|\tilde \phi_0''(\t)|^2}{\tilde
   \phi_0(\t)}\, {\mathrm d}\t < \infty.
   \end{matrix}
  $$

Finally, we obtain the estimate
 \beq
 \label{886}
  \mbox{$
  \frac 52 \int\limits_0^T \int\limits_0^L u^2 \phi_0 \, {\mathrm d}x\, {\mathrm d}t \le
  \frac{c_0L^7}{7T^3} \, - J_0, \quad \mbox{where}
  \quad J_0= \int\limits_0^Lu_t(x,0)(L-x)^3 \, {\mathrm d}x.
   $}
   \eeq
 Thus, if $J_0 >0$,  the solution must blow-up in finite time $T_0$
satisfying
 $
  \mbox{$
  T_0 \le \big(\frac{c_0L^7}{7 J_0}\big)^{\frac 13}.
   $}
   $
It is curious that this blow-up conclusion depends on second data
$u_t(x,0)$ and does not involve the first initial function
$u(x,0)$.

\section{Shock and rarefaction similarity solutions for $S_\mp$, $H(\pm)$, or others}
 \label{Sect2}

As a natural next step, before proposing concepts on existence,
uniqueness, and entropy description of solutions, one needs to get
a detailed understanding of the types of singularities that can be
generated by NDEs. As often happens in nonlinear evolution PDEs,
the refined structure of such bounded or unbounded shocks
 can be described by similarity solutions, which thus  we have to begin with.

\subsection{Finite time blow-up formation of the shock wave
$S_-(x)$}
 \label{Sect3.1}

\begin{figure}
\centering
\includegraphics[scale=0.75]{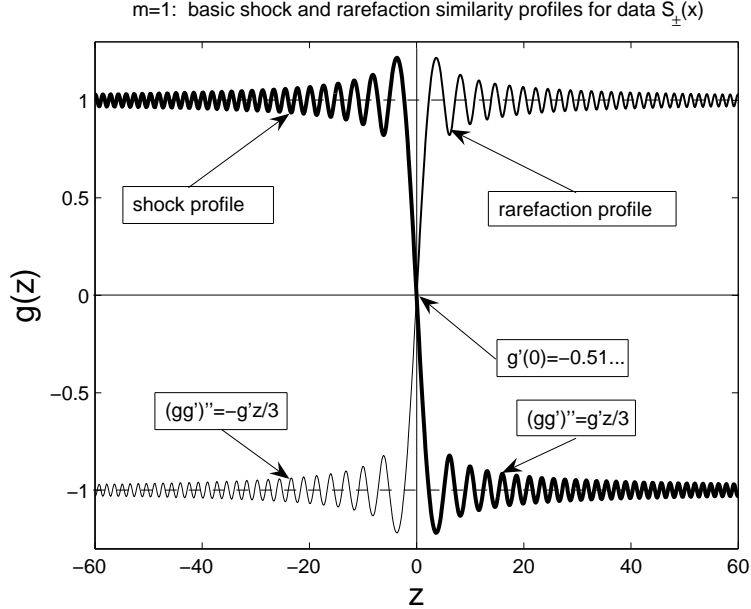}
 \vskip -.3cm
\caption{\small The shock similarity profile as the unique
solution of the problem (\ref{2.2}).}
   \vskip -.3cm
\label{F1}
\end{figure}

To this end, we use  the following similarity solution of the NDE
(\ref{1}):
 \beq
 \label{2.1}
 u_-(x,t)=g(z), \quad z= x/(-t)^{\frac 13},
  \eeq
where $g$ solves the ODE problem
 \beq
 \label{2.2}
  \mbox{$
  (g g')''= \frac 13 \, g'z \quad \mbox{in} \quad \re, \quad f(\mp
  \infty)=\pm 1.
   $}
   \eeq
   By translation, the blow-up time in (\ref{2.1}) reduces to $T=0$,
so that we want, in the sense of distributions or in
 $L^1_{\rm loc}$ (see Proposition \ref{Pr.Co} below),
 \beq
 \label{2.3}
  u_-(x,t) \to S_-(x) \quad \mbox{as}
  \quad t \to 0^-.
   \eeq
    In view of the symmetry of the ODE (\ref{2.2}),
     \beq
     \label{symm88}
      \left\{
      \begin{matrix}
     g \mapsto -g,\\
     z \mapsto -z,
      \end{matrix}
      \right.
      \eeq
it suffices to study the odd solutions for $z<0$ with the
anti-symmetry conditions,
 \beq
 \label{2.4}
 g(0)=g''(0)=0.
  \eeq
A typical structure of this shock similarity profile  $g(z)$
solving (\ref{2.2}) is shown in Figure \ref{F1} by the bold line.

\ssk

\noi\underline{\em On regularization in numerical methods}. Since
we are going  to use essentially the numerical results now and
later on,
especially for
higher-order NDEs, let us
 describe the peculiarities of such numerics.
 For the third-order equations such as (\ref{2.2}),
this and further numerical constructions are performed by {\tt
MatLab} with the standard {\tt ode45} solver therein.
 Currently, we use the enhanced
relative and absolute tolerances
 \beq
 \label{T1NN}
 {\rm Tols}=10^{-12}.
  \eeq
Instead of the degenerate ODE (\ref{2.2}), we have solved the
regularized equation (this idea is naturally associated with the
general concepts
of parabolic $\e$-approximation 
\cite[\S~7]{GalEng})
 \beq
 \label{T2}
  \mbox{$
  g''' =\frac {{\rm sign}\, g}{\sqrt{\nu^2+g^2}} \, \big(\frac 13 \, g'z- 3 g'
  g''\big),
 \quad \mbox{with the regularization parameter $\nu=10^{-12}$},
 $}
  \eeq
  where the choice of small $\nu$ is coherent with the
  tolerances in (\ref{T1NN}).

\ssk

\noi\underline{\em On asymptotics, existence, and uniqueness}. One
can see that $g(z)$ is always oscillatory about constant
equilibria $\mp 1$ as $z \to \pm \infty$. Indeed, for $z \ll -1$,
where $g(z) \approx 1$, the ODE (\ref{2.2}) asymptotically reduces
to the linear equation
 \beq
 \label{2.5}
  \mbox{$
  g'''= \frac 13 \, g'z + ... \, .
   $}
   \eeq
Hence,  $g(z)$ satisfies the asymptotics of
 the classic Airy function:
 as $z \to
 -\infty$,
 \beq
 \label{Ai.1}
  \mbox{$
g(z) \sim 1+ c {\rm Ai}(z) \sim 1+c  |z|^{-\frac 14} \cos\big(a_0
|z|^{\frac 32}+c_0\big), \,\,\,\, \mbox{where} \,\,\,\, a_0= \frac
29 \,\sqrt
 3.
  $}
 \eeq
Concerning the original PDE (\ref{1}), this means that, as $x \to
-\infty$, and hence $u \to 1$, the NDE is  asymptotically
  transformed into the {\em linear dispersion equation}
 \beq
 \label{2.6}
 u_t=u_{xxx}, \quad \mbox{with the fundamental solution}
  \eeq
  \beq
  \label{bb123}
  \mbox{$
b(x,t)= t^{-\frac 13} {F}\big( x/{t^{\frac 13}}\big),
\quad\mbox{where} \,\,\, F={\rm Ai}(z), \quad F'' + \frac 13 \,F
\,z=0, \quad \int F=1.
 $}
 \eeq

In particular,
 the asymptotics (\ref{Ai.1}) implies that
 the {\em total
 variation} (TV) of such solutions of (\ref{1}) (and hence $u_-(x,t)$ for any $t<0$)
  is {\em infinite}. Setting $|z|^{3/2}=v$ in the integral
  yields
 \beq
 \label{pp1ss}
 \mbox{$
 | g(\cdot)|_{\rm TV} =
 \int\limits_{-\infty}^{+\infty} |g'(z)|\, {\mathrm d}z \sim \int\limits^\infty \frac
 {|\cos v|}{v^{1/6}}\, {\mathrm d}v= \infty.
 $}
  \eeq
 This is in striking contrast with the case of conservation laws
 (\ref{3}), where finite total variation approaches and  Helly's second theorem
 (compact embedding of sets of bounded functions of bounded total variations into $L^\infty$)
   are  key, \cite{Ol1}.
  In view of the
  presented  properties of the similarity profile $g(z)$, the
  convergence in (\ref{2.3}) takes place
 for any $ x \in \re$, uniformly in $\re \setminus (\d,\d)$, $\d>0$ small, and
 in $L^p_{\rm loc}(\re)$ for  $p \in[1, \infty)$; see below.


Before passing to more accurate mathematical treatment of the
profile $g(z)$, we present its  regular asymptotic expansion near
the origin: for any $C <0$, there exists a unique solution of the
ODE (\ref{2.2}), (\ref{2.4}) satisfying
 \beq
 \label{2.7}
  \mbox{$
 g(z) = C z+ \frac 1{72}\, z^3+
  \frac 1{72^2}\, \frac 1C \, z^5 +... \quad
 (C=-0.51... \,\,\,\mbox{for $S_-$; see Figure
 \ref{F1}}).
  $}
  \eeq
The uniqueness of such  asymptotics is traced out by using
Banach's  Contraction Principle
 applied to the equivalent integral equation in the metric
 of $C(-\d,\d)$, with $\d>0$  small.

 In addition, we will use the following scaling invariance of the
ODE in (\ref{2.2}): if $g_1(z)$ is a solution, then
 \beq
 \label{2.8}
 \mbox{$
 g_a(z) = a^3 g_1\big(\frac z a\big) \quad \mbox{is a solution for any $a \not =
 0$}.
  $}
  \eeq

\begin{proposition}
\label{Pr.1} The problem  $(\ref{2.2})$   admits the unique shock
wave profile $g(z)$, which is an odd analytic function and is
strictly positive for $z<0$.
 \end{proposition}

\begin{figure}
\centering
\includegraphics[scale=0.65]{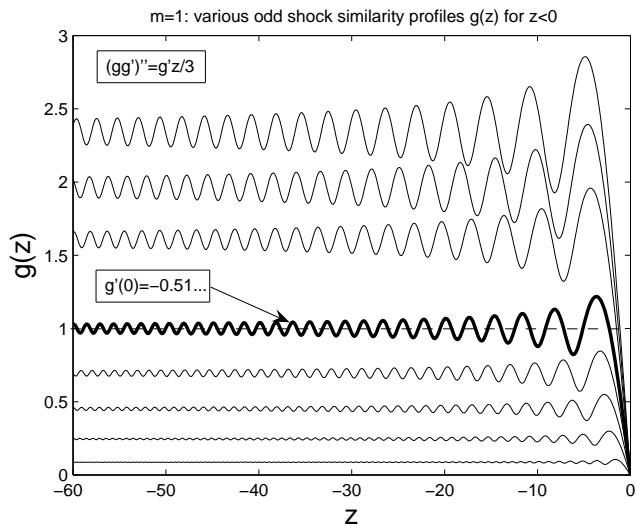}
 \vskip -.3cm
\caption{\small The odd shock similarity profiles $g(z)$ for $z<0$
with various limits $g(z) \to C_-$ as $z \to - \infty$.}
   \vskip -.3cm
\label{F2}
\end{figure}

Uniqueness follows from the asymptotics (\ref{2.7}) and scaling
invariance (\ref{2.8}). Global existence as infinite extension of
the unique solution from $z=0^-$ follows from the ODE (\ref{2.2}),
which, besides solutions with finite limits as $z \to -\infty$,
admits the unbounded solutions with the behaviour
 \beq
 \label{2.9}
  \mbox{$
  g(z)= \frac 1{60} \, z^3 +... \quad \mbox{for} \quad |z| \gg 1.
 $}
   \eeq
 Figure \ref{F2} shows other
profiles from this one-parameter family with different limits of
$g(z)$ as $z \to -\infty$. Of course, all of them are related to
each other by scaling (\ref{2.8}).

The proof of the positivity is rather technical and will be
focused on later. In particular, our reliable numerics with
enhanced accuracy show the exhaustive positivity of the similarity
profile $g(z)$, and, actually,  this can be treated as a
computational-based proof.

Thus, we need to prove analyticity at $z=0$ only, since by
positivity $g(z)$ is analytic at any other point. Differentiating
(\ref{2.2}) $k-2$ times yields by Leibnitz formula that the
derivatives at $z=0$ satisfy the recursion relation
 $$
  \mbox{$
  (k+1)g' g^{(k)}= 2(k-2)g^{(k-2)} - \sum\limits_{i=2}^{k+1}
  g^{(i)}g^{(k+1-i)} \frac{(k+1)!}{i! \, (k+1-i)!}.
 $}
 $$
It is not difficult to conclude that the growth rate of the
derivatives
 $
  |g^{(k)}(0)|$ is not essentially more than $ k!$.
This means that the corresponding power series
 $$
 \mbox{$
 g(z) = \sum\limits_{(k)} \frac {g^{(k)}(0)}{k!} \, z^k
  $}
  $$
  has the unit radius of convergence.

\ssk

We now return to the convergence (\ref{2.3}).

\begin{proposition}
 \label{Pr.Co}
 For the shock similarity profile $g(z)$ from
 Proposition $\ref{Pr.1}$,
   $(\ref{2.3})$:

 {\rm (i)} does not hold in $L^1(\re)$, and

  {\rm (ii)} does hold in $L^1_{\rm loc}(\re)$.

  \end{proposition}

  \noi{\em Proof.} (i) follows from (\ref{pp1ss}), since
   $\, g(z)-1 \not \in L^1(\re_-)$.
 (ii) Here we need an extra estimate by using the asymptotics (\ref{Ai.1}):
  for a fixed finite $l>0$, as $t \to 0^-$,
  \beq
  \label{ee1ss}
  \begin{matrix}
  \int\limits_{-l}^0|g(z)-1| \, {\mathrm d}x=
   (-t)^{1/3}\int\limits_{-l(-t)^{-1/3}}^0|g(z)-1| \, {\mathrm d}z
   \ssk\ssk\\
 \sim (-t)^{1/3}\int\limits_{-l(-t)^{-1/3}}^0 |z|^{-\frac 14}|\cos(a_0 |z|^{\frac 32})|
  \, {\mathrm d}z
 \ssk\ssk\\
\sim (-t)^{1/3}\int\limits_0^{l^{3/2}(-t)^{-1/2}} v^{-\frac
12}|\cos(a_0 v)| \, {\mathrm d}v \sim (-t)^{\frac 1{12}} \to 0.
\quad \qed
  \end{matrix}
  \eeq

\ssk

Thus, the rate of convergence in (\ref{ee1ss}) and the fact of
convergence itself depend on the delicate asymptotics of the Airy
function, and more precisely, on the tail structure of the
fundamental solution of the corresponding linearized equation
(\ref{2.6}). Therefore, for NDEs such as (\ref{1}) or other
higher-order and non-fully divergent ones, the topology of
convergence {\em cannot be obtained in a unified manner} and is
individual. For instance, this can be $L^p_{\rm loc}$-convergence,
where $p>1$ may depend on the NDE under consideration.
Fortunately, due to (\ref{ee1ss}) (note that the ``gap" for such
convergence expressed by the positive exponent $\frac 1{12}$ is
rather small), for (\ref{1}), we still have the
 $L^1_{\rm loc}$-convergence as for the conservation law
 (\ref{3}), where this topology is naturally reinforced by Helly's
 second theorem on compact embedding of BV into $L^\infty$.


\subsection{On formation of non-symmetric final time profile}

\begin{figure}
\centering
\includegraphics[scale=0.65]{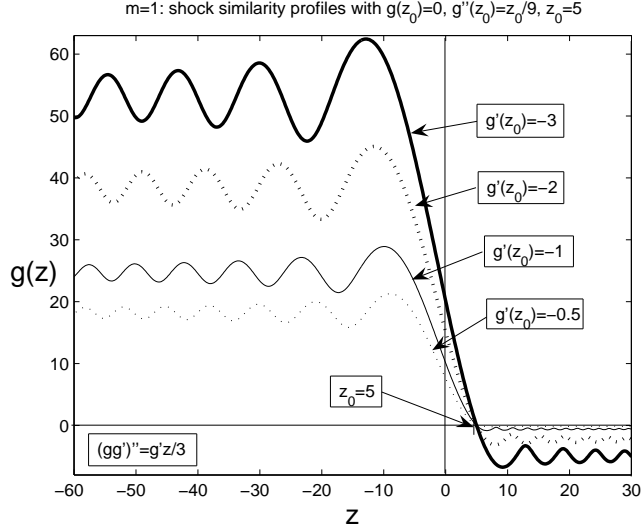}
 \vskip -.3cm
\caption{\small Non-symmetric shock similarity profiles $g(z)$
with various limits $g(z) \to C_\pm$ as $z \to \pm \infty$.}
   \vskip -.3cm
\label{F3}
\end{figure}

This analysis is harder but repeats the above arguments. Figure
\ref{F3} shows a
few of such similarity profiles $g(z)$ that,
according to (\ref{2.1}), generate non-symmetric step-like
functions, so that, as $t \to 0^-$,
 \beq
 \label{2.10}
  u_-(x,t) \to
  \left\{
  \begin{matrix}
  C_->0 \,\,\, \mbox{for} \,\,\, x<0, \\
  C_0  \,\,\,\quad \, \quad \mbox{for} \,\,\, x=0, \\
 C_+<0 \,\,\, \mbox{for} \,\,\, x>0,
  \end{matrix}
   \right.
    \eeq
where $C_- \not = -C_+$ and $C_0 \not = 0$. In order to understand
the whole variety of such non-symmetric profiles, one needs to
check
 how many regular orbits can pass through any singular point
 $z=z_0>0$, at which $g(z_0)=0$ ($z_0=5$ in Figure \ref{F3}). By Banach's Contraction
 Principle, it can be shown that, for any $z_0>0$, there exists a
 1D family (a bundle), with the regular expansion
  \beq
  \label{2.11}
 \mbox{$
  g(z)=C(z-z_0) + \frac {z_0}{18} \, (z-z_0)^2 + \frac 1{72} \,
  (z-z_0)^3+... \, , \quad \mbox{where $C <0$}.
   $}
    \eeq
For $z_0=0$, this coincides with (\ref{2.7}).

Thus, the total family of profiles regularly (moreover,
analytically) passing through singular points $z=z_0$ is 2D, with
parameters $z_0$ and $C$. By scaling invariance (\ref{2.8}), this
variety can be reduced to a 1D manifold with, say, fixed $z_0=1$
and arbitrary first derivative $C<0$.

\subsection{Shock similarity profiles with finite interfaces: weak discontinuities}

These correspond to $C_+=0$ in (\ref{2.10}). Especially, we are
interested in checking how to get in the blow-up limit as $t \to
0^-$ the reflected Heaviside function (\ref{2.12}). The value at
the origin
    $H(0)\in (0,1)$ is obtained via integration.
  It follows from (\ref{2.11}) with $C=0$ that such
regular weak solutions of the equation in (\ref{2.2}) have the
following asymptotics:
 \beq
 \label{2.13}
  \mbox{$
 g(z)= \frac {z_0}{18} \, (z-z_0)^2+... \,\,\,\mbox{as} \,\,\, z
 \to z_0^-,
  $}
  \eeq
where $z_0>0$ is the only free parameter. Proof of existence of
such profiles, which are positive for all $z<z_0$, is similar to
that of Proposition \ref{Pr.1}.

\begin{figure}
\centering
\includegraphics[scale=0.65]{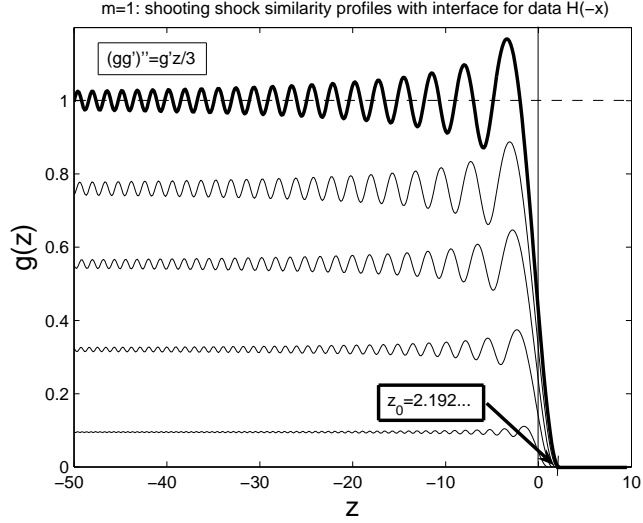}
 \vskip -.3cm
\caption{\small The shock similarity profiles $g(z)$ with finite
interfaces at $z=z_0>0$.}
   \vskip -.3cm
\label{F4}
\end{figure}

 Using the
additional smoothness of expansion (\ref{2.13}), we set $g(z)
\equiv 0$ for $z > z_0$, thus creating an admissible type of
proper {\em weak discontinuity} at $z=z_0$ for the ODE (and the
NDE), where the flux $(g g')'(z)$ is continuous and $g'(z)$ is
Lipschitz.
 Several profiles
from this family are shown in Figure \ref{F4}. The bold line
corresponds to formation of the Heaviside shock (\ref{2.12}).
 A mathematical justification of the
consistency of this  cutting-off procedure (actually meaning
finite propagation of the right-hand interface) will be performed
later. Here we note the following approximation property:
 the
profiles with the finite interface at a singular point $z=z_0$
 \beq
 \label{ss0}
  \begin{matrix}
\mbox{can be obtained as the limit as $C_+ \to 0^-$}\\
 \mbox{of smooth
 solutions (as in Figure $\ref{F3}$).}
 \end{matrix}
  \eeq
Note that the family with $C_+<0$ are composed from analytic
solutions, though these have singular points where $g=0$. On the
other hand, we expect that $g(z)$ with interface at $z_0>0$ can be
obtained as the limit $\mu \to 0^+$ of strictly positive solutions
$\{g_\mu(z)\}$, where $g_\mu(z) \to \mu>0$ as $z \to +\infty$.
This demands extra asymptotic analysis of matching the flow
(\ref{2.2}) about $g=+1$ as $z \to -\infty$ with the 1D stable
manifold of $g=\mu$ as $z \to +\infty$.

\begin{figure}
\centering
\includegraphics[scale=0.65]{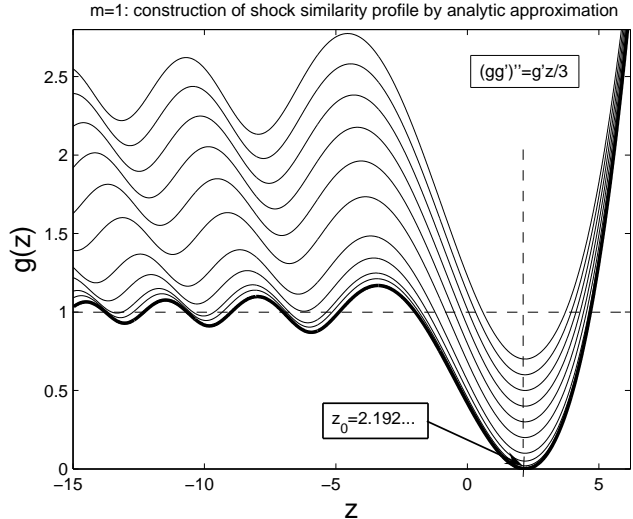}
 \vskip -.3cm
\caption{\small Analytic approximation of $g(z)$ with the
behaviour (\ref{2.13}) at $z=z_0>0$.}
   \vskip -.3cm
\label{F41}
\end{figure}

 As an alternative, Figure
\ref{F41} shows another way of approximation of  $g(z)$ with zero
at $z=z_0$ by a family of strictly positive analytic solutions of
the ODE (\ref{2.2}). This convinces that the solutions with the
singular point at $z=z_0$ and behaviour (\ref{2.13}) (but not with
finite interface!) are G-admissible in the ODE sense; see Section
\ref{Sect5}.

\subsection{Rarefaction similarity solutions}

Using the reflection symmetry
 \beq
 \label{symm1}
 \left\{
 \begin{matrix}
u \mapsto -u, \\
 t \mapsto -t, \,
  \end{matrix}
  \right.
  \eeq
   of the PDE (\ref{1}), we have
that it admits  similarity solutions defined for all $t>0$,
 \beq
  \label{2.14}
  u_+(x,t) = g(z), \,\,\, \mbox{with} \,\, z=x/t^{\frac 13},
   \eeq
where now $g(z)$ solves
 \beq
 \label{2.2R}
  \mbox{$
  (g g')''=- \frac 13 \, g'z \quad \mbox{in} \quad \re, \quad f(\mp
  \infty)=\mp 1.
   $}
   \eeq
Obviously, these profiles are obtained from the blow-up ones in
(\ref{2.2}) by reflection,
 \beq
 \label{2.15}
\mbox{if shock profile $g(z)$ solves (\ref{2.2}), then $g(-z)$ is
rarefaction one in (\ref{2.2R}).}
 \eeq
The corresponding rarefaction similarity profile is shown in
Figure \ref{F1}. This regular similarity solution of (\ref{2.2R})
has  the necessary initial data: by Proposition \ref{Pr.Co}(ii),
in $L^1_{\rm loc}$,
 \beq
 \label{2.3R}
  u_+(x,t) \to S_+(x) \quad \mbox{as}
  \quad t \to 0^+.
   \eeq
Other profiles $g(-z)$ from shock wave similarity patterns
generate further rarefaction solutions including those with finite
left-hand interfaces.

Thus, all three Riemann's problems for the NDE (\ref{1}) with data
(\ref{Ri1}) and (\ref{2.12})
 admit a unified
  similarity treatment.

\section{Unbounded shocks and other self-similar singularities}
 \label{Sect3}

 We continue to introduce other types of shocks and
 singularities that are associated with the NDE (\ref{1}).

\subsection{Blow-up self-similar solutions: invariant subspace and
 critical exponent $\a_{\rm c}=-\frac 1{10}$}

We now consider  more general blow-up similarity solutions of
 (\ref{1}),
 \beq
 \label{3.1}
  \mbox{$
 u_\a(x,t)=(-t)^\a g(z), \quad z= x/(-t)^{\b}, \quad \b=
 \frac{1+\a}3 \quad(\a \in \re),
  $}
  \eeq
where $g$ solves the ODE
 \beq
 \label{3.2}
  \mbox{$
 {\bf A}(g) \equiv (g g')''= \frac {1+\a}3 \, g'z -\a g \equiv {\bf C}g \quad \mbox{in} \quad \re.
   $}
   \eeq
In (\ref{3.1}), we introduce an extra arbitrary parameter $\a \in
\re$. We next show that the behaviour of the similarity profiles
$g(z)$ and hence of the corresponding solutions $u_\a(x,t)$
essentially depend on whether $\a>0$ or $\a<0$.

We first prove the following auxiliary result  explaining a key
feature of the ODE (\ref{3.2}).

\begin{proposition}
 \label{Pr.2}
{\rm (i)} Both nonlinear ${\bf A}$ and linear ${\bf C}$ operators
in
 $(\ref{3.2})$ admit the $4D$ linear space
  \beq
  \label{3.3}
  W_4={\rm Span}\{1,z,z^2,z^3\}.
   \eeq

   {\rm (ii)} The ODE $(\ref{3.2})$ possesses  nontrivial solutions on $W_4$
   in  two cases:
 \beq
 \label{3.4}
  \mbox{$
 {\rm (I)} \quad
 \fbox{$
 \a=\a_{\rm c}=- \frac 1{10},
 $}
 \quad \mbox{with the solutions  given by}
 $}
  \eeq
 \beq
 \label{3.5}
  \mbox{$
 g(z)= C_0 + C_1z + \frac 1{60} \, z^3, \quad \mbox{where $C_{0,1} \in
 \re$ are arbitrary constants, and}
 $}
  \eeq
   \beq
 \label{3.41}
  \mbox{$
 {\rm (II)} \quad\a =-1, \quad \mbox{with}
 \quad
 g(z)= \frac{400}3\, C_2^3 + 20 C_2^2 z + C_2 z^2+ \frac 1{60} \,
 z^3,\,\,\, C_2 \in \re.
 $}
  \eeq
 \end{proposition}

 \noi{\em Proof.} (i) is straightforward, since, for any
  \beq
  \label{3.6}
  g=C_0 + C_1z + C_2z^2+ C_3 z^3 \in W_4,
   \eeq
   the following holds:
    \beq
    \label{3.7}
     \begin{matrix}
    {\bf A}(g)= 6(C_1C_2+C_0C_3) + 12 (C_2^2+2 C_1C_3)z+ 60 C_2C_3
    z^2+ 60 C_3^2 z^3 \in W_4,\qquad \ssk\ssk\\
    {\bf C}g= - \a C_0+ \frac{1-2\a}3\, C_1 z + \frac{2-\a}3\,
    C_2z^2+ C_3 x^3 \in W_4. \qquad\qquad \qquad\qquad\qquad
    \qquad
     \end{matrix}
     \eeq

(ii) According to the ODE (\ref{3.2}), equating the coefficients
given in (\ref{3.7}) yields the algebraic system
 \beq
 \label{3.8}
 \left\{
 \begin{matrix}
6(C_1C_2+C_0C_3)=-\a C_0,\,\,\,\\
  12 (C_2^2+2 C_1C_3)=\frac{1-2\a}3 \,C_1,\\
60 C_2 C_3  =\frac{2-\a}3 \, C_2,\,\,\,\,\,\qquad\quad\\
 60 C_3^2=C_3.\qquad\qquad\qquad\,\,\,\,
  \end{matrix}
  \right.
   \eeq
The last equation yields $C_3= \frac 1{60}$ (we exclude the easy
case $C_3=0$), and then the third one implies that either $C_2  =
0$, or
 \beq
 \label{3.9}
 \mbox{$
60 C_2 \frac 1{60} = \frac{2-\a}3 \, C_2 \quad \Longrightarrow
\quad \a=-1.
 $}
 \eeq
Assuming first that $C_2=0$ and substituting into the second
equation in (\ref{3.8}), we infer that, for arbitrary $C_1 \not =
0$,
 \beq
 \label{3.10}
  \mbox{$
 24 C_1 \frac 1{60} =\frac{1-2\a}3 \,C_1
 \quad \Longrightarrow \frac 25=\frac{1-2\a}3, \,\,\,
 \mbox{i.e.,\,\,
 $\a=-\frac 1{10}$}.
 $}
 \eeq
The first equation is also valid for any $C_0 \in \re$.

Choosing $C_2 \not = 0$ leads to the less interesting case
(\ref{3.41}).
 $\qed$

 \ssk

 \noi{\bf Remark: blow-up for the NDE on the invariant subspace.}
 Since the linear subspace
 (\ref{3.3}),  $W_4={\rm Span}\{1,x,x^2,x^3\}$
  is invariant under the quadratic operator ${\bf A}$, the whole original NDE
 (\ref{1}) can be restricted to $W_4$, on which the evolution of
 the corresponding solutions,
  \beq
  \label{sol1}
  u(x,t)= C_0(t) + C_1(t) x+ C_2(t)x^2 + C_3(t) x^3 \in W_4,
  \eeq
  is governed by a dynamical system with the
  right-hand sides from (\ref{3.8}),
 \beq
 \label{3.8NN}
 \left\{
 \begin{matrix}
C_0'=6(C_1C_2+C_0C_3),\\
  C_1'=  12 (C_2^2+2
C_1C_3),\\
 C_2'=60 C_2 C_3,\,\,\,\,\,\qquad\quad\\
C_3'= 60 C_3^2.\qquad\qquad\,\,\,\,\,
  \end{matrix}
  \right.
   \eeq
This is easily integrated starting from bottom and gives the
blow-up behaviour,
 $$
 \begin{matrix}
 C_3(t)= \frac 1{60(T-t)}, \,\,\, C_2(t)= \frac {A_0}{T-t},
  \,\,\, C_1(t)=\frac {B_0}{(T-t)^{2/5}}-\frac {60}{(T-t)^{3/5}}, \ssk\ssk\\
  C_0(t)= \frac {D_0}{(T-t)^{1/{10}}}+\frac {20 A_0 B_0}{(T-t)^{
  2/{5}}}
  - \frac {720 A_0}{(T-t)^{3/{5}}},
   \end{matrix}
$$
 where the blow-up time $T$,
 $A_0$, $B_0$, and $C_0$ are arbitrary constants.
Unlike the above gradient blow-up, the present singularity
formation is governed by the cubic growth of the solution
(\ref{sol1}) as $x \to \infty$, so that it is less interesting.

\subsection{Blow-up similarity solutions and their
properties for $\a \in (\a_{\rm c},0)$}
 \label{S3.2}

\begin{figure}
\centering
\includegraphics[scale=0.75]{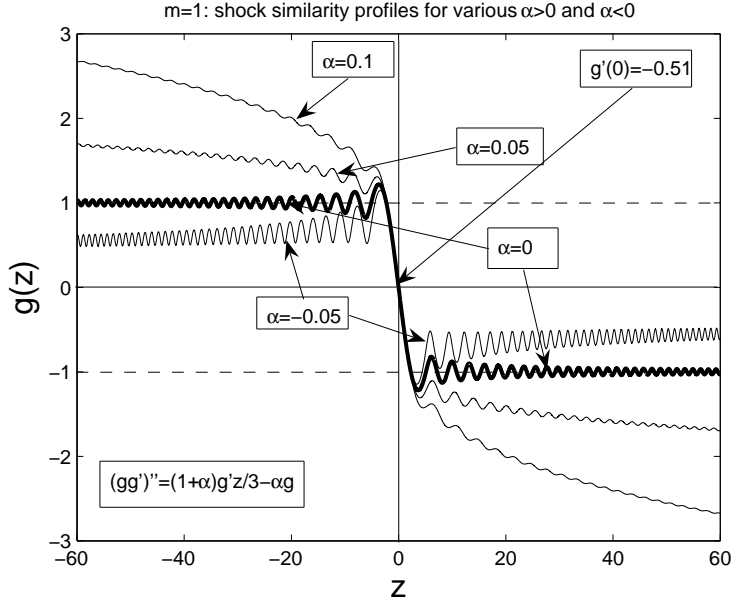}
 \vskip -.3cm
\caption{\small The shock similarity profiles as  solutions of the
ODE (\ref{3.2}) for $\a>0$ and $\a<0$.}
   \vskip -.3cm
\label{F5}
\end{figure}

Basic local mathematical properties of the ODE (\ref{3.2}) are
similar to those for (\ref{2.2}), so we omit these details and
concentrate on principal global features of blow-up patterns. The
typical structure of solutions $g(z)$ of (\ref{3.2}) can be
understood from Figure \ref{F5}, where we present a few profiles
for $\a>0$ and $\a<0$, including the bold line for $\a=0$ that has
been studied in the previous section. All these profiles satisfy
the anti-symmetry conditions (\ref{2.4}).

For $\a>0$, the ODE (\ref{3.2}) suggests (and Figure \ref{F5}
confirms this) that the similarity profiles are growing as $z \to
-\infty$ according to the linear part of the equation, i.e.,
 \beq
 \label{3.11}
  \mbox{$
  \frac {1+\a}3 \, g'z -\a g \approx  0\,\,\, \mbox{for
$z \ll -1$}  \,\,\, \Longrightarrow
  \,\,\, g(z)= C_-|z|^{\frac{3\a}{1+\a}} + ... \, , \,\,\, C_->0.
  $}
   \eeq
Passing to the limit $t \to 0^-$ in (\ref{3.1}) does not lead to a
shock wave but  to a ``weak singularity",
 \beq
 \label{3.12}
 u_\a(x,0^-)=  C_\mp|x|^{\frac{3\a}{1+\a}} \quad \mbox{for $x<0$ and
 $x>0$}.
 \eeq
For $\a= \frac 12$ or $\a=2$, these give the behaviour $\sim x$ or
$\sim |x|x$ near the origin respectively.

\ssk

The behaviour for $\a<0$ is more interesting, and deserves further
analysis. First of all, for $\a \in (-\frac 1{10},0)$, the
similarity profiles $g(z)$ are still positive for $z<0$ and have
the asymptotic behaviour (\ref{3.11}). Therefore, (\ref{3.12}) is
also true but one should take into account that here the exponent
 $\frac{3\a}{1+\a}<0$, so that the final time profile is
 unbounded, i.e.,
  \beq
  \label{3.91}
 u_\a(x,0^-)=  C_\mp|x|^{\frac{3\a}{1+\a}} \to \pm \infty \quad \mbox{as $x \to 0^\mp$}.
 \eeq
This is the first example of strongly discontinuous {\em
unbounded} shocks that can be obtained evolutionary via blow-up of
sufficiently smooth continuous self-similar solutions of the NDE
(\ref{1}). The conservation law (\ref{3}) cannot produce such a
shock, since by the Maximum Principle,
 $$
 \mbox{$
  \sup_x |u(x,t)| \le \sup_x |u_0(x)| \quad \mbox{for all} \quad t>0.
   $}
   $$

\begin{figure}
\centering
\includegraphics[scale=0.75]{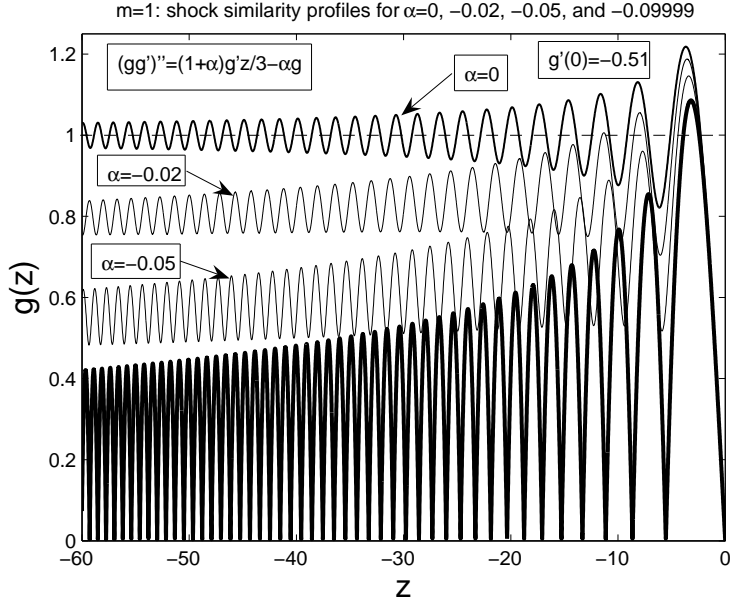}
 \vskip -.3cm
\caption{\small Shock similarity profiles as  solutions of the ODE
(\ref{3.2}) for some $\a \in (-\frac 1{10},0]$; the {\bf boldface
line} corresponds to $\a=-0.09999$ ($\approx \a_{\rm c}=-\frac
1{10}$ for the ``saw" profile).}
   \vskip -.3cm
\label{F6}
\end{figure}

\begin{figure}
\centering
\subfigure[convergence to $g_{\rm c}(z)$]{
\includegraphics[scale=0.52]{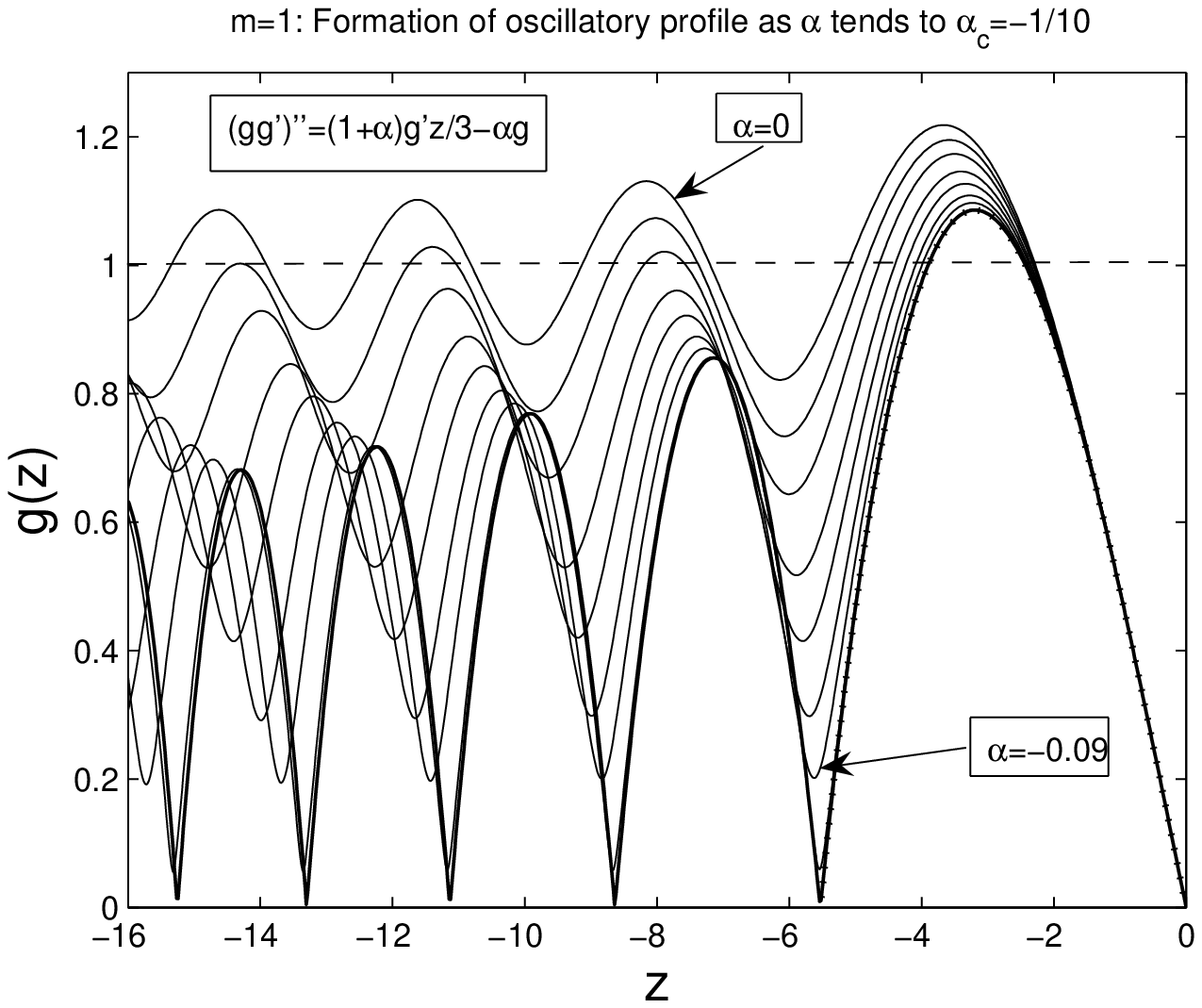} 
}
\subfigure[zeros of $g_{\rm c}(z)$]{
\includegraphics[scale=0.52]{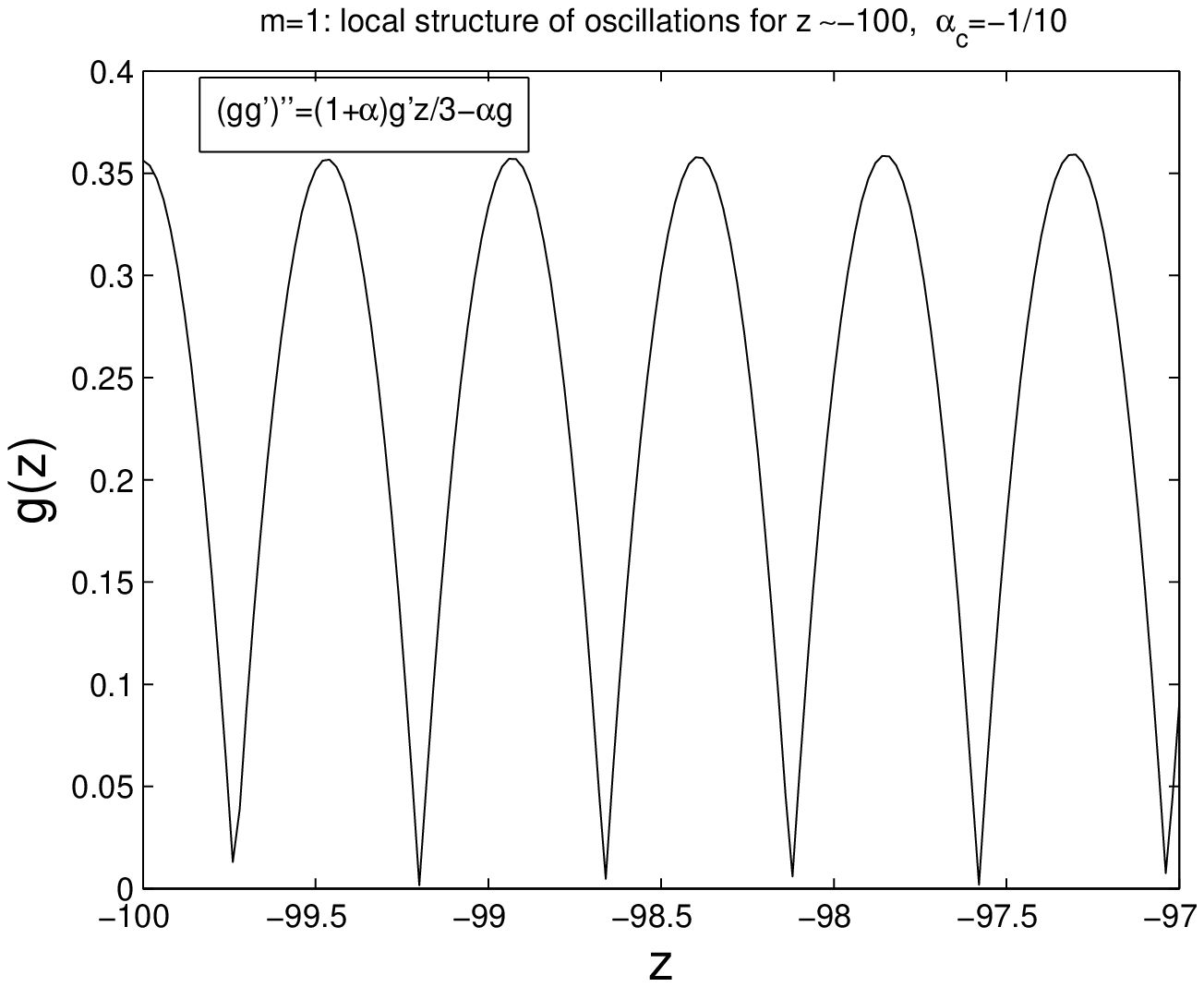} 
}
 \vskip -.4cm
\caption{\rm\small The limit (\ref{3.13}); convergence to the
singular profile $g_{\rm c}(z)$ (a), and  the asymptotic structure
of the ``tail" of $g_{\rm c}(z)$ for $z \sim -100$.}
 \vskip -.3cm
 \label{F7}
\end{figure}

 In Figure \ref{F6}, we show how the positive similarity
profile $g(z)$ for $z < 0$ is deformed according to the limit
 \beq
 \label{3.13}
  \mbox{$
g(z) \to g_{\rm c}(z) \quad \mbox{as} \quad \a \to \a_c^-=- \frac
1{10}.
  $}
  \eeq
 Note that a slightly easier mechanism of formation of such ``saw cusps" is a well-known
bifurcation phenomenon in the FFCH  (\ref{R1}) and
Degasperis--Procesi (\ref{DP0}) type equations; see
\cite[p.~422]{Shen06} and references therein.  Figure \ref{F7}(a)
explains the actual convergence of positive profiles in the limit
(\ref{3.13}), while (b) shows the oscillatory  structure of
$g_{\rm c}(z)$ for $z \ll -1$ at $\a=-\frac 1{10}$.

\begin{figure}
\centering
\subfigure[$C_- \not =-C_+$]{
\includegraphics[scale=0.52]{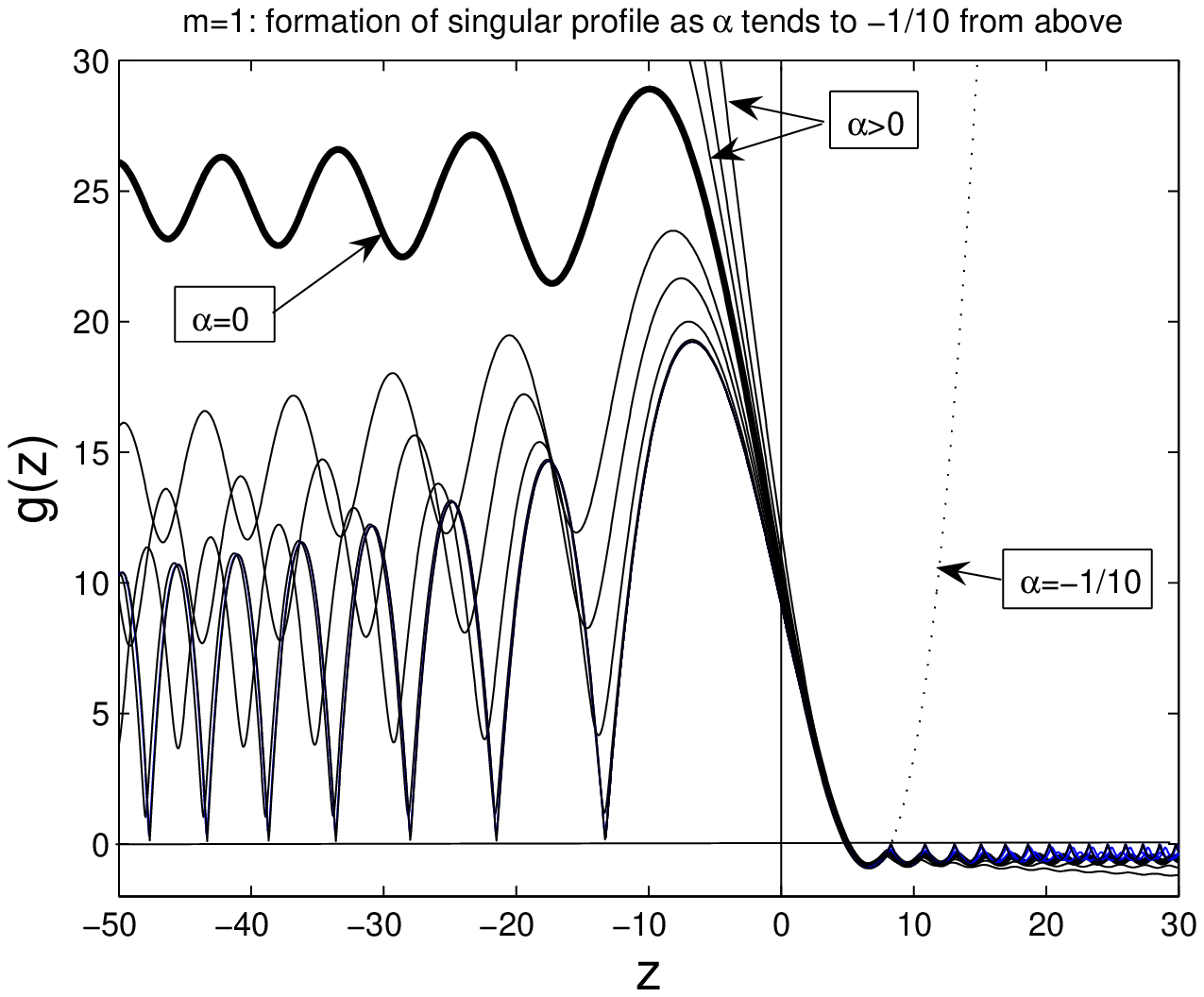} 
}
\subfigure[with finite interface]{
\includegraphics[scale=0.52]{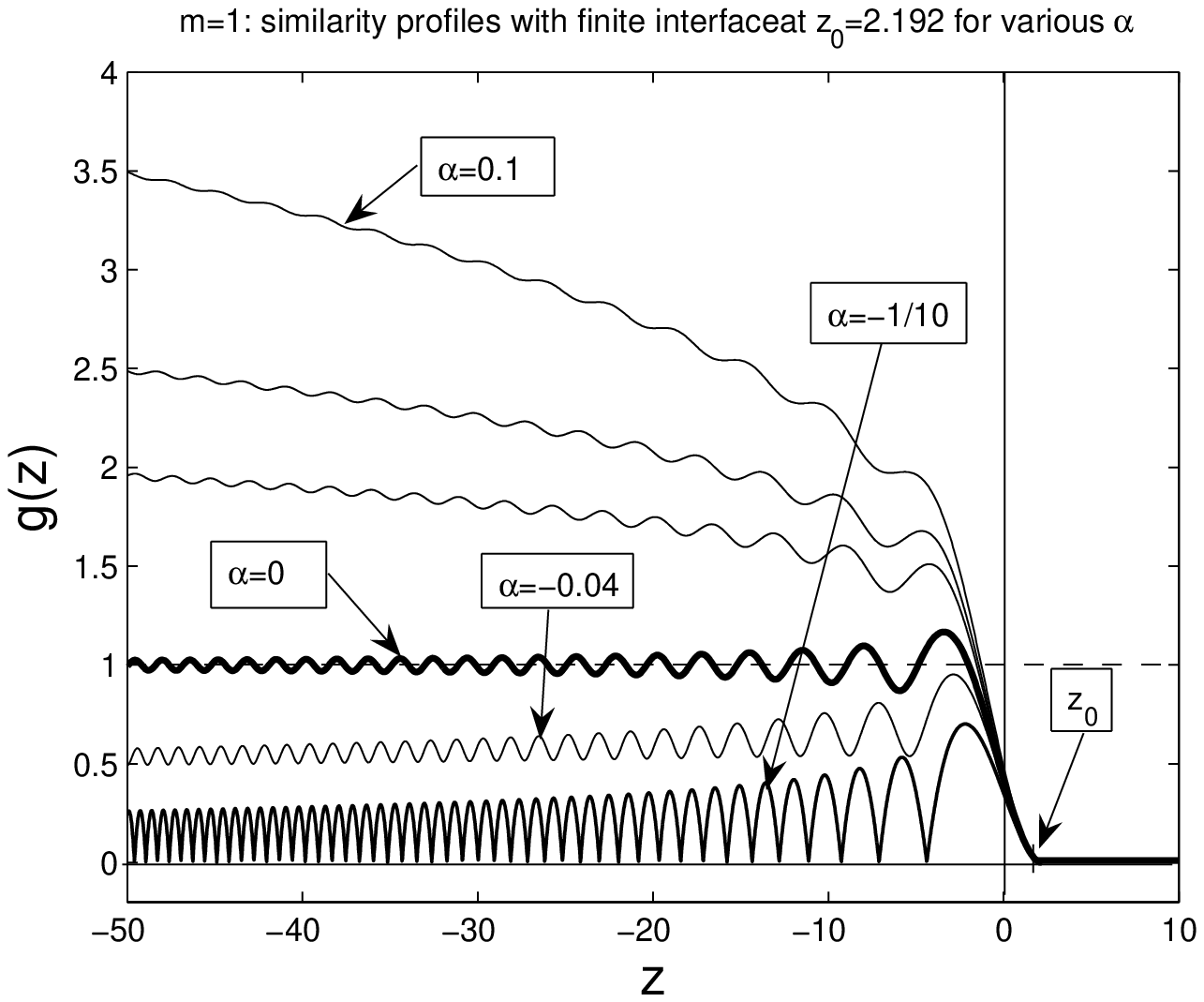} 
}
 \vskip -.4cm
\caption{\rm\small Various similarity solutions of (\ref{3.2}) for
positive and negative $\a$.}
 \vskip -.3cm
 \label{F55}
\end{figure}

Finally, in Figure \ref{F55}, we show other examples of
non-symmetric profiles $g(z)$ for various $\a$ including those in
(b) that  have the finite interface fixed  at $z_0=2.192$. The
bold profile generates as $t \to 0^-$ the Heaviside function
 (\ref{2.12}), with parameters
 \beq
 \label{pp1}
 H(0)=0.4197..., \quad z_0=2.192... \, .
  \eeq


\ssk

\noi{\bf Remark: complete blow-up.} There is an important question
concerning continuation of the similarity solution satisfying
(\ref{3.91}) beyond blow-up, i.e., for $t>0$. This is not the
subject of the
 study, so we present a few comments on that.
According to extended semigroup theory of blow-up solutions (see
\cite{GalGeom} for references and known results for
reaction-diffusion PDEs), it is natural to construct such an
extension via a smooth approximation (truncation) of the NDE by
introducing a bounded nonlinearity,
 \beq
 \label{nn1}
  \mbox{$
 u_n: \quad u_t =( \varphi_n(u))_{xxx}, \quad \mbox{with}
 \quad \varphi_n(u) = \frac 12\, \frac{n^2 u^2}{n^2+u^2} \to \frac 12 \,
 u^2
 \,\,\, \mbox{as} \,\,\, n \to \infty.
  $}
  \eeq
Taking the same initial data we obtain a sequence $\{u_n(x,t)\}$
of better bounded solutions (bounded since blow-up demands the
quadratic growth $\sim u^2$ as $u \to \iy$ of the coefficient that
is not available anymore; we omit difficult technical details).

Nevertheless, each $u_n(x,t)$ can produce a shock at some $t=t_n
\approx 0$, but now of finite height $\sim 2n$.  Since
$u_{nt}(x,t_n)$ gets very large in a neighbourhood of $x=0^-$
(this follows from the fact that, for the original solution,
$u_t(x,0^-)\to +\infty$ as $x \to 0^-$ for ``data" (\ref{3.91})),
one can expect that $u_n(x,t)$ for $x<0$ grows  rapidly for small
$t>0$ achieving fast the value $\sim n$. This means that
$u_n(x,t)$ for small $t>0$ gets the shape of $n S_-(x)$.

Finally, we then conclude that for any fixed $t>0$, as $n \to
\infty$,
 \beq
 \label{nn2}
 u_n(x,t) \sim n S_-(x) \to \left\{
 \begin{matrix}
 +\infty \quad \mbox{for} \,\,\, x<0, \ssk\\
 -\infty \quad \mbox{for} \,\,\, x>0.
  \end{matrix}
  \right.
   \eeq
Actually, this means that this function $\bar u(x,t) = \pm \infty$
is the correct proper solution as a natural continuation of the
similarity solution (\ref{3.1}), (\ref{3.91}) for $t>0$. In other
words, we observe {\em complete blow-up} at $t=0^-$ and the fact
that the extended semigroup of proper (minimal) solutions is
discontinuous at this blow-up time.

\subsection{Oscillatory invariant profile for $\a_{\rm c}=- \frac 1{10}$}

Then the ODE (\ref{3.2}) takes the form
 \beq
 \label{3.99}
  \mbox{$
 (g g')''= \frac 3{10} \, g'z+ \frac 1{10} \, g,
  $}
  \eeq
which is associated with the invariant subspace $W_4$
 in (\ref{3.3}).

\begin{figure}
\centering
\subfigure[two profiles]{
\includegraphics[scale=0.52]{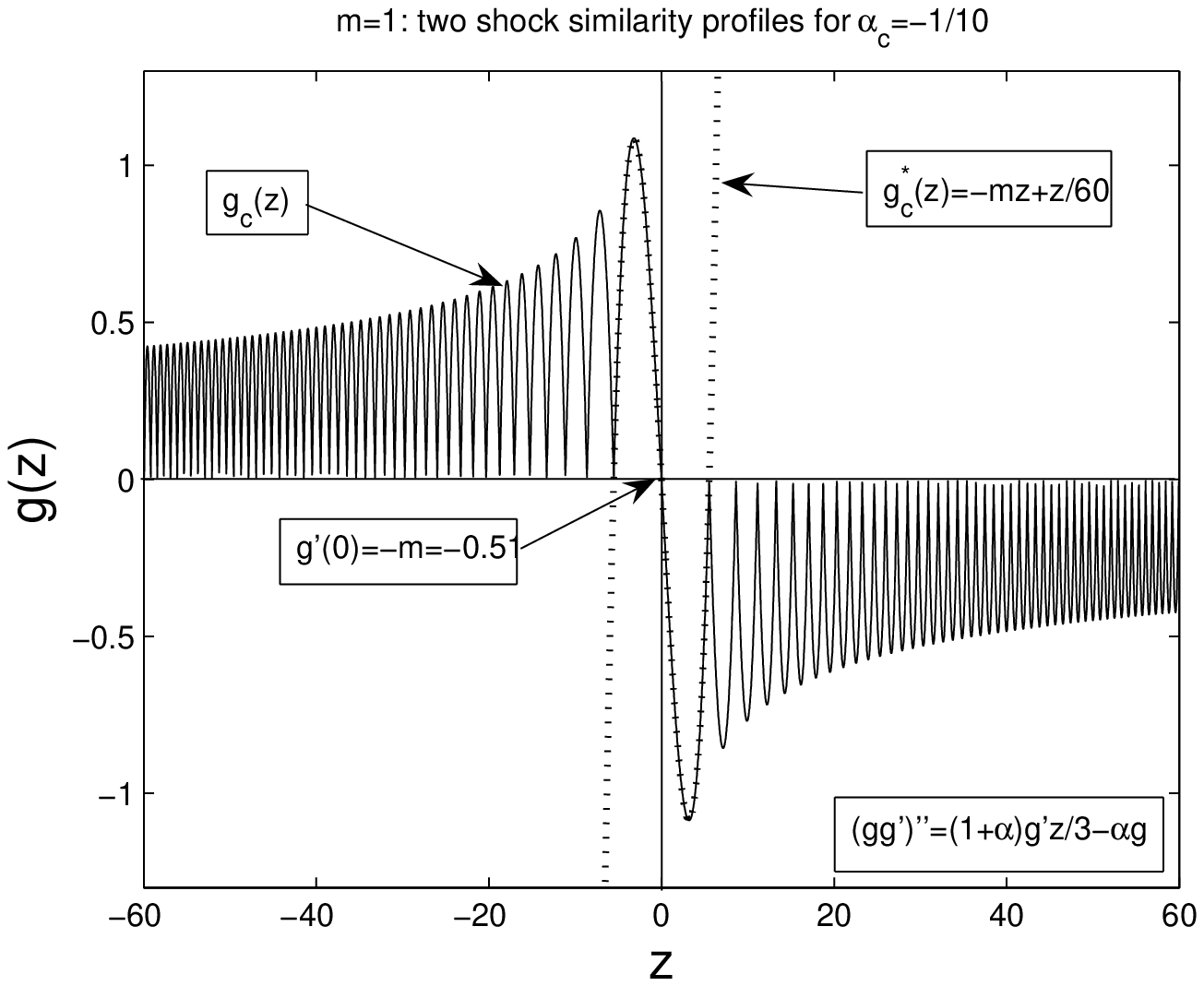} 
}
\subfigure[envelope to the saw]{
\includegraphics[scale=0.52]{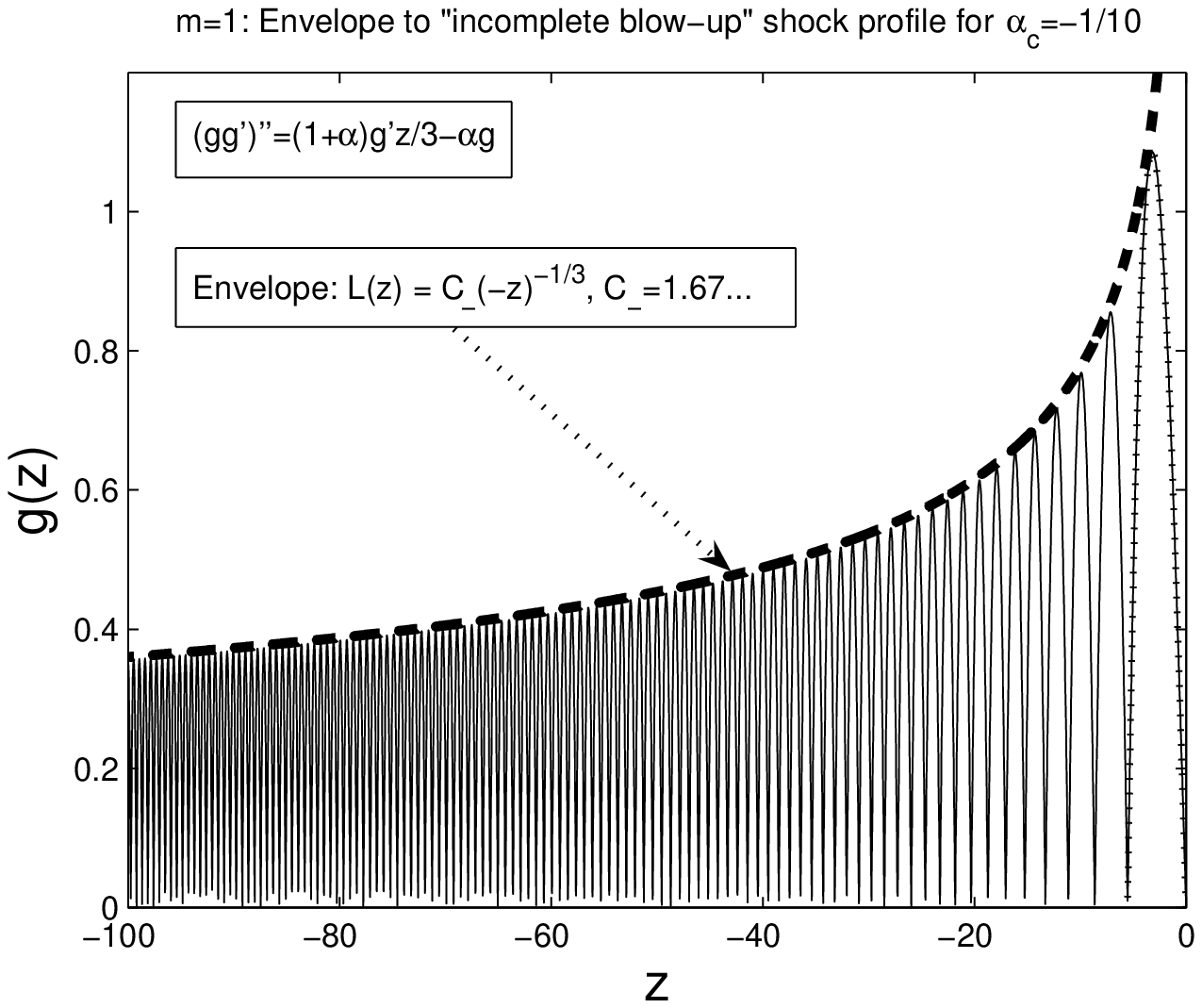} 
}
 \vskip -.4cm
\caption{\rm\small Two solutions of (\ref{3.99}) on $W_4$.}
 \vskip -.3cm
 \label{F8}
\end{figure}

 We have already seen, that the
critical case $\a=\a_{\rm c}$ is of special interest.  Figure
\ref{F8}(a) shows that, in fact, for $\a_{\rm c}=- \frac 1{10}$,
there exist two continuous solutions. The second one given by the
dotted line is the exact invariant one
 \beq
 \label{3.14}
  \mbox{$
  g_{\rm c}^*(z) = -m z + \frac 1{60} \, z^3\in W_4, \quad \mbox{where}
  \quad  -m= g'(0)<0 .
  $}
  \eeq
Figure \ref{F8}(b) explains the structure of the {\em envelope} to
the main ``saw-type" oscillatory profile $g_{\rm c}(z)$ that has
the same structure as in (\ref{3.11}), i.e.,
 \beq
 \label{3.01}
 L(z)=C_- |z|^{-\frac 13} \quad \mbox{as}
 \quad z \to -\infty.
  \eeq
Therefore, the final time profile $u_\a(x,0^-)$ is not like
(\ref{3.91}). Indeed, the envelope $L(x)$ only has this form, so
that we have another topology of convergence:
 \beq
 \label{3.02}
  u_\a(x,t) \rightharpoonup C_- L(x) \quad \mbox{as} \quad t \to 0^-
   \eeq
   in the weak sense in $L^\infty$ or in $L^2$. Here $C_-=1.67...$ is a
   constant that characterizes the oscillatory part of this shock
   wave
   distribution; see more details below.


\begin{figure}
\centering
\includegraphics[scale=0.65]{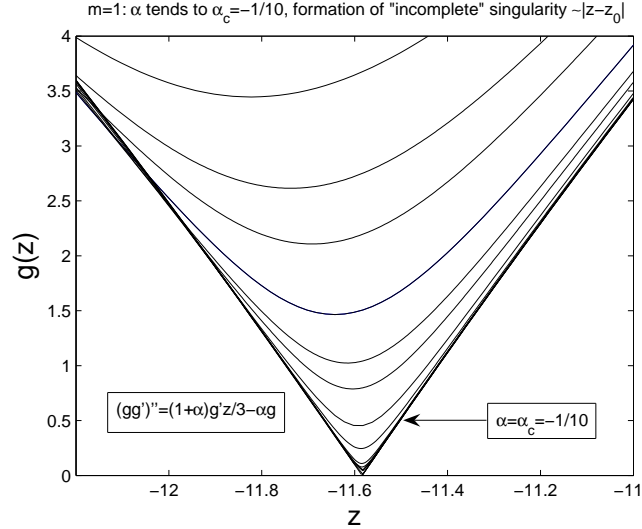}
 \vskip -.3cm
\caption{\small Local convergence as $\a \to -\frac 1{10}$ of
smooth solutions $g(z)$ to the ``saw" $g_{\rm c}(z)$ with the
behaviour (\ref{3.51}).}
   \vskip -.3cm
\label{F9}
\end{figure}

Figure \ref{F9} explains  formation of  zero at some $z=z_0<0$ of
$g_{\rm c}(z)$ as $\a \to -\frac 1{10}$. It is clearly seen that
$g_{\rm c}(z)$ has a symmetric behaviour for $z \approx z_0$, so
instead of (\ref{2.11}), 
 \beq
 \label{3.51}
  \mbox{$
g_{\rm c}(z)= C|z-z_0|+ \frac{z_0}{18}\,(z-z_0)^2+... \, , \quad
C>0.
$}
 \eeq
 Since the ``flux" $(g_{\rm c}^2)''$ is continuous at $z=z_0$,
 (\ref{3.51}) represents true  weak solutions of (\ref{3.99}).
 Moreover,
it is also crucial for our further analysis that
 \beq
 \label{ss1}
 \begin{matrix}
 \mbox{the weak
piece-wise smooth solution $g_{\rm c}(z)$ can be obtained}\\
 \mbox{ as the
limit of $C^3$
solutions of the family of ODEs (\ref{3.2})}.
 \end{matrix}
  \eeq

Finally, we prove the following result concerning the invariant
properties of the oscillatory profile $g_{\rm c}(z)$:

\begin{proposition}
\label{Pr.4} For $\a=- \frac 1{10}$,  in the piece-wise sense,
 \beq
 \label{3.52}
 g_{\rm c}(z) \in W_4.
  \eeq
  \end{proposition}

  \noi{\em Proof.} Without loss of generality, we consider the
  solution $g_{\rm c}(z)$ of (\ref{3.99}) with $m=1$, 
 \beq
 \label{3.53}
 g(0)=g''(0)=0 \quad \mbox{and} \quad g'(0)=-1.
  \eeq
  By uniqueness, it follows from (\ref{3.5}) that the first hump
  of the profile shown in Figure \ref{F7}(a) is explicitly given
  by
   \beq
   \label{3.54}
    \mbox{$
   g_{\rm c}(z)=-z + \frac 1{60} \, z^3 \quad \mbox{for} \quad  z \in (z_0,0),
   \quad \mbox{where} \quad  z_0=-\sqrt{60} \quad(g_c(z_0)=0).
    $}
    \eeq
For the second hump with the transition given by (\ref{3.51}), we
take the ``full" profile (\ref{3.5}) and demand two matching
conditions at $z=z_0$,
 $$
 g(z_0)=0, \quad g'(z_0^+)=2=-g'(z_0^-), \quad \mbox{i.e.,}
 $$
 \beq
 \label{3.55}
 \left\{
 \begin{matrix}
 C_0+C_1 z_0+ \frac 1{60} \, z_0^3=0, \ssk\\
 C_1+ \frac 1{20} \, z_0^2=-2.\qquad\,
  \end{matrix}
  \right.
   \eeq
   This yields
    $
   C_1=-5$
   and
   $C_0=-4 \sqrt{60},
    $
    so that the second hump is given by
     \beq
     \label{3.56}
 \mbox{$
     g_{\rm c}(z)= - 4 \sqrt{60}-5z+ \frac 1{60} \, z^3
     \quad \mbox{on} \quad (z_1,z_0).
     $}
     \eeq
The second zero $z_1$ is
 \beq
 \label{3.57}
  \mbox{$
 z_1=- \rho \sqrt{60} \,\, \Longrightarrow \,\, \rho^3- 5 \rho
 +4=0, \,\,\, \mbox{i.e.,} \,\,\, \rho=
 \frac{\sqrt{17}-1}2=1.56155... \,.
  $}
  \eeq
This  procedure is continued without bound and gives the ``saw" as
in Figure \ref{F8}(b). $\qed$

\subsection{Local complete blow-up in the ODE for $\a<\a_{\rm c}$}

\begin{figure}
\centering
\includegraphics[scale=0.75]{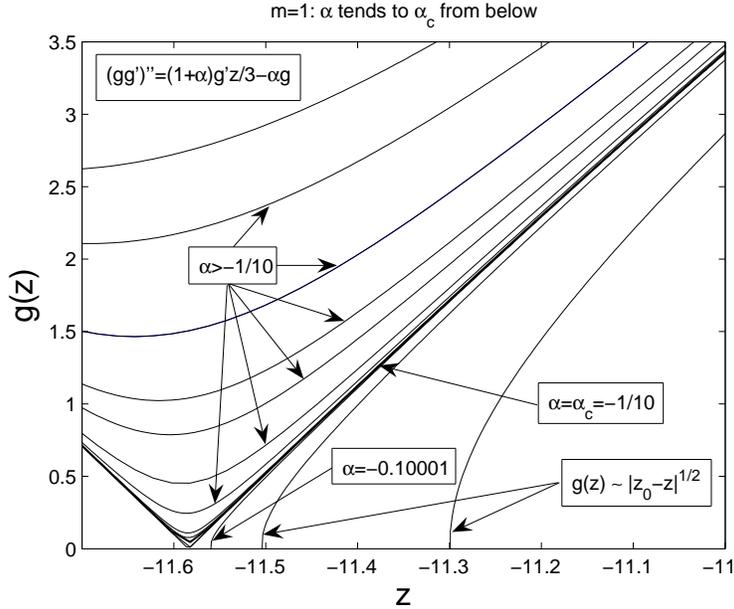}
 \vskip -.3cm
\caption{\small Complete blow-up of orbits for $\a < -\frac
1{10}$.}
   \vskip -.3cm
\label{F10}
\end{figure}

It turns out that for $\a < - \frac 1{10}$, the orbits $g(z)$
cannot pass through singularity at $\{g=0\}$. This local
``complete blow-up" is clearly illustrated by the equation for the
second derivative
 \beq
 \label{3.58}
  \mbox{$
 v=g'' \quad \Longrightarrow \quad g v'' +  v^2 + 4 g' v'=
 \frac{1+\a}3 \, v z + \frac{1-2\a}3\, g'.
 $}
  \eeq
In other words, this blow-up $v \to -\infty$ as $z \to z_0^+<0$ is
performed according to the quadratic source term,
 \beq
 \label{v1}
  \mbox{$
 v'' =- \frac 1g \, v^2+... \, .
  $}
  \eeq
The ODE (\ref{3.2}) for $\a \not = -1$ admits the following
behaviour at such complete blow-up zeros:
 \beq
 \label{v2}
 \mbox{$
 g(z)=C\sqrt{z-z_0}+\frac {4(1+\a)}{45} \, z_0(z-z_0)^2+... \quad \mbox{as} \quad z \to z_0^+, \quad
 C>0.
  $}
  \eeq
  For $\a=-1$, the remainder is $\sim O((z-z_0)^3)$.
  The results of such  shootings  are presented in Figure
  \ref{F10}.

\subsection{Rarefaction similarity solutions}


These are defined as follows:
 \beq
  \label{2.14RR}
  u_\a(x,t) = t^\a g(z), \,\,\, \mbox{with} \,\,\,\, z=x/t^{\frac {1+\a}3},
  \quad \mbox{where}
   \eeq
 \beq
 \label{2.2RR}
  \mbox{$
  (g g')''=- \frac {1+\a}3 \, g'z +\a g \quad \mbox{in} \quad \re.
   $}
   \eeq
This gives the reflected shock blow-up profiles $g(-z)$, and all
the above properties can be translated to the rarefaction
evolution.

In connection with this, we mention that the unbounded initial
data (see (\ref{3.91}))
 \beq
 \label{in1}
 u_0(x)= \left\{
 \begin{matrix}\,\,
 x^{\frac{3\a}{1+\a}} \qquad \mbox{for} \,\,\, x>0,\\
-|x|^{\frac{3\a}{1+\a}} \quad \mbox{for} \,\,\, x<0,
 \end{matrix}
 \right.
  \eeq
  for any $\a \in (-\frac 1{10},0)$ generate a bounded (for any $t>0$) similarity
  solution (\ref{2.14RR}) with the corresponding profile $g(-z)$.
  According to our further theory, this
  is an example of an entropy solution that is expected to be unique.

\section{On G-admissible solutions of ODEs}
 \label{Sect5}

The approximating  properties that we underlined in (\ref{ss0})
and (\ref{ss1}) are key for understanding the ODEs involved.
Actually, these are related to
 the concept of
 {\em G-admissibility} of solutions of ODEs with shocks. It was introduced in 1959 by Gel'fand
\cite{Gel} and was developed on the basis of TW-solutions of
hyperbolic equations and systems; see details in \cite[\S~2,
8]{Gel}.

Namely, talking for definiteness  the typical ODEs (\ref{2.2}) or
(\ref{3.2}), a non-classical (i.e., not $C^3$) solution is called
G-admissible if it can be constructed by a converging sequence of
smooth solutions of the same ODE. In particular, this means that
the profiles with finite interfaces as in Figure \ref{F4} and even
the ``saw-type" profile $g_{\rm c}$ for $\a_{\rm c}=- \frac 1{10}$
are G-admissible.

Of course, this is an ODE concept, but we will bear it in mind
 when develop a PDE concept of entropy solutions in
 \cite{GalNDE5, 3NDEII}.

\section{On travelling wave and generic formation of moving shocks}
 \label{Sect31}

We now briefly discuss how to construct moving shock waves, since
previously we concentrated on standing shocks such as $S_-(x)$,
$H(-x)$, and others. Of course, the first simple idea is to
consider:

\subsection{Travelling waves (TWs).} These are solutions of (\ref{1}) of
the form
 \beq
 \label{T1}
 u(x,t)=f(y), \quad y=x-\l t,
  \eeq
  where $\l \in \re$ is the TW speed. Substituting into (\ref{1})
  and integrating once yields
   \beq
   \label{T2N}
   - \l f'=(ff')'' \quad \Longrightarrow \quad - \l f= (ff')'+A_0,
    \eeq
    where $A_0 \in \re$ is an arbitrary constant.
If a shock (a finite discontinuity) occurs at some point $y_0 \in
\re$, (\ref{T2N}) yields the following {\em Rankine--Hugoniot}
 {\em condition} for the speed:
 \beq
 \label{T20}
  \mbox{$
 \l=- \frac{[(ff')'](y_0)}{[f](y_0)}\, ,
 $}
  \eeq
where, as usual, $[(\cdot)](y_0)$ stands for the jump of the
function $(\cdot)(y)$ at $y_0$.
In the case
 \beq
 \label{T3}
 \l=A_0=0,
 \eeq
 we obtain either constant solutions
  \beq
  \label{T4}
  f(y)=C, \quad \mbox{and hence $f=S_\pm(y)$, $H(\pm y)$, and
  others, \,\,\, or}
   \eeq
 \beq
 \label{T5}
  f(y)= \sqrt {A_1 y + A_2}, \quad \mbox{where} \quad A_{1,2} \in
  \re,
   \eeq
 from which one can reconstruct various discontinuous solutions satisfying
the Rankine--Hugoniot condition (\ref{T20}). At $y_0= -\frac
{A_2}{A_1}$ (cf. (\ref{v2})), the solution (\ref{T5}) is not
sufficiently smooth to be treated as a solutions of the Cauchy
problem.

For $\l \not = 0$, (\ref{T2N}) admits smoother explicit solutions
\beq
 \label{T6}
  \mbox{$
  f(y)= -\frac \l 6 \, (y+B)^2 - \frac 3{2\l}\, A_0, \quad B \in \re,
  $}
  \eeq
which also give various shocks. In general, it is easy to describe
all the solutions admitted by the standard second-order equation
(\ref{T2N}) that we are not doing here and will focus on the
following principle question:

\subsection{On formal theory of generic formation of finite moving shocks}
 \label{Sect5.2}

Let a uniformly bounded continuous  weak solution $u(x,t)$ of the
NDE (\ref{1}) in $\re \times (-1,0)$ (for convenience, the
focusing blow-up time is again $T=0$) create a discontinuous shock
at
 \beq
 \label{T61}
 x=0, \,\,\, t=0^-, \quad \mbox{at which the shock moves with the speed
 $\l$}.
  \eeq
Then it is natural to consider this phenomenon within the moving
frame and  introduce the corresponding TW variable $y=x - \l t$ as
in (\ref{T1}). Hence,   the solution $u=u(y,t)$ now solves
 \beq
 \label{T7}
 u_t = (uu_y)_{yy} + \l u_y \quad \mbox{for} \quad y \in \re,
 \,\,\, t \in(-1,0).
  \eeq

We now study formation of bounded shocks  associated with the
similarity variable $z$  in (\ref{2.1}). Indeed, there are other
 blow-up
scenarios, with variables as in (\ref{3.1}); see \cite{Gal3NDENew}
for  more generic single point ``gradient catastrophes" (these
require other rescalings). Hence, according to (\ref{2.1}), we
introduce the rescaled variable
 \beq
 \label{T8}
 u(y,t)= v(z, \t), \quad \mbox{where}
 \quad z= x/(-t)^{\frac 13}, \,\,\, \t = - \ln (-t) \to
 +\infty\,\, \mbox{as}
 \,\, t \to 0^-.
  \eeq
  Substituting this into (\ref{T7}) yields that
 $v(z,\t)$ solves the following non-autonomous PDE:
  \beq
  \label{T9}
   \mbox{$
  v_\t = (v v_z)_{zz} - \frac 13 \, v_z z + \l{\mathrm e}^{-\frac 23
  \, \t} g_z \quad \mbox{in} \quad \re \times \re_+.
   $}
   \eeq
Here the right-hand side is an exponentially small for $\t \to
\infty$ perturbation of the stationary operator of the ODE in
(\ref{2.2}). Therefore, it is natural to expect that the
asymptotic behaviour as $\t \to +\infty$ (i.e., $t \to 0^-$) is
described by the non-perturbed equation,
 \beq
 \label{T10}
  \mbox{$
 v_\t = (v v_z)_{zz} - \frac 13 \, v_z z= {\bf B}(v) \equiv {\bf A}(v)- {\bf
 C}v \quad \mbox{in} \quad \re \times \re_+.
  $}
  \eeq

Therefore, stabilization in (\ref{T9}) and (\ref{T10})  to a
self-similar profile,
 \beq
 \label{T11}
 v(z,\t) \to g(z) \quad \mbox{as \,$\t \to
\infty$ uniformly on compact subsets},
 \eeq
 is equivalent to the  question of the stable manifold  of the stationary solution $g(z)$.
  In the linear approximation, this
 will depend on the spectrum of the non self-adjoint third-order
 operator
  \beq
  \label{T12}
   \mbox{$
  {\bf B}'(g)Y= (g Y)''' - \frac 13 \, Y'z.
 $}
   \eeq
This spectral problem is not easy, though similar higher-order
linear operators already occurred in some evolution odd-order
PDEs; see \cite[\S~5.2]{GalEng},  \cite[\S~9.2]{2mSturm}
for details.

Considering, for simplicity, ${\bf B}'(g)$ in the space of odd
functions defined in $\re_-=\{z<0\}$, in view of conditions
(\ref{2.4}) on $g(z)$ and the smooth behaviour (\ref{2.7}), the
end point $z=0$ can be treated as a regular one. At infinity,
where $g(-\infty)=1$, we arrive asymptotically at the linear
operator as in (\ref{2.5}), which admits Airy function as one of
its eigenfunctions. These operators are naturally defined in the
weighted space
 \beq
 \label{T13}
 L^2_\rho(\re_+), \quad \mbox{where} \quad
 \rho(z)={\mathrm e}^{-a|z|^{3/2}}
  \eeq
and $a>0$ is a sufficiently small constant; see more details in
\cite[\S~9.2]{2mSturm}.

It is important that, in view of the scaling  symmetry
(\ref{2.8}),
 \beq
 \label{T14}
  \begin{matrix}
 \l=0 \in \s({\bf B}'(g)) \quad
\mbox{with the eigenfunction} \ssk\ssk\\
 \psi_0(z)= \frac{\mathrm
d}{{\mathrm d}a} \, g_a(z)\big|_{a=1}= 3 g(z)- g'(z) z \in
L^2_\rho,
\end{matrix}
 \eeq
 so that ${\bf B}'(g)$ has the non-empty centre subspace $E^{\rm
 c}={\rm Span}\{\psi_0\}$.

 Then, according to typical trends of invariant manifold theory; see \cite{Lun}
 (we must admit that sectorial properties
 of ${\bf B}'(g)$ in $L^2_\rho$ are still unknown, while compactness
  of the resolvent is most plausible),  it is natural first to
 try a
 centre subspace behaviour for (\ref{T10}).
Stability on the centre manifold is  key for understanding the
stabilization properties of the flow (\ref{T10}). To this end, we
perform the ``linearization" procedure
 \beq
 \label{T15}
v(z,\t)=g(z)+Y(z,\t)
 \eeq
  to get the following equation with a quadratic perturbation:
   \beq
 \label{T151}
 \mbox{$
Y_\t= {\bf B}'(g)Y + \frac 12 \, (Y^2)_{zzz}.
 $}
 \eeq
In order to formally check the stability of the centre manifold,
we assume that for $\t \gg 1$ the behaviour follows the centre
subspace, i.e.,
  \beq
 \label{T16}
 \mbox{$
Y(z,\t)= a(\t) \psi_0(z) + w(z,\t), \quad \mbox{where}
 \quad w \bot E^{\rm c}.
 $}
 \eeq
Here, the orthogonality is defined according to the corresponding
Eucledian  space with the indefinite metric $\langle \cdot,\cdot
\rangle_*$,
 $$
 \mbox{$
\langle v,w \rangle_*= \int v(z) w(-z)\, {\mathrm d}z;
 $}
 $$
 see details and examples in \cite[\S~9.2]{2mSturm}. Substituting
(\ref{T16}) into (\ref{T151}) and multiplying by the corresponding
normalized adjoint eigenfunction $\psi_0^*$ yields
 \beq
 \label{T17}
 \mbox{$
a'=-\g_0 a^2+... \,\,\, \mbox{as} \,\,\, \t \to \infty, \quad
\mbox{where} \quad \g_0= \frac 12 \,\langle (\psi_0^2)''',\psi_0^*
\rangle_*,
 $}
 \eeq
 where the scalar products involved are defined
 by the extension of the linear functionals by
 Hahn--Banach's classic theorem (we omit the details that are not crucial
 here).
For existence of such a centre subspace behaviour, it is necessary
that
 \beq
 \label{T18}
  \mbox{$
\g_0 \not =0 \quad \Longrightarrow \quad \exists \,\,\, a(\t)=
\frac 1{\g_0 \, \t} + ... \to 0 \quad \mbox{as} \quad \t \to
\infty.
 $}
 \eeq


On the other hand, stabilization in (\ref{T11}) is also possible
along a stable infinite-dimensional subspace of $\BB'(g)$. Note
that this would be the only case if $\g_0=0$ in (\ref{T18}). Such
operators (without degeneracy at $z=0$, which nevertheless is not
strong and/or essential) are known to admit a point spectrum
\cite[\S~9.2]{2mSturm}. Such formal stable manifold behaviours
admit similar standard calculus, though the duality between the
operator  $\BB'(g)$ and the adjoint one $\BB'^*(g)$, as usual for
odd-order singular operators, will provide us with several
technical difficulties. Recall again that there exist other more
structurally stable (and possibly generic) ways to create single
point gradient blow-up singularities, where other rescaled
nonlinear and linearized operators occur, \cite{Gal3NDENew}. We do
not touch those questions here.


Overall, this explains  some formalities of the asymptotic
stability analysis, which remains a difficult open problem. Even checking the necessary non-orthogonality assumptions,
such as (\ref{T18}), can be a difficult problem, both analytically
and  numerically.

Anyway, it is key that if  stability of stationary profiles $g(z)$
for (\ref{T10}) and (\ref{T9}) takes place (at least, on some
manifold) then, for
 bounded moving shocks, we can use a similar geometric rule
(\ref{rr1})  distinguishing $\d$-entropy  and non-entropy
discontinuous solutions (see more  details in \cite{3NDEII}).



\enddocument